\documentclass[dvips]{amsart}
\usepackage{math,pb-diagram}
\newcommand\smashnode[2][1]{\node[#1]{\phantom{[}\quad\makebox[0mm]{$#2$}\quad\phantom{]}}}

\begin{document}
\title{On Parabolic Subgroups and Hecke Algebras of Some Fractal Groups}
\date{May 11, 2000; revised February 22, 2001}
\author{Laurent Bartholdi}
\email{laurent@math.berkeley.edu}
\address{11, ch. de la Barillette, 1260 Nyon, Switzerland}
\author{Rostislav I. Grigorchuk}
\email{grigorch@mi.ras.ru}
\address{Steklov Mathematical Institute, Gubkina 8, 117966 Moscow,
  Russia}
\thanks{The second author wishes to express his thanks to the ``Swiss
  National Science Foundation''}
\begin{abstract}
  We study the subgroup structure, Hecke algebras, quasi-regular
  representations, and asymptotic properties of some fractal groups of
  branch type.

  We introduce parabolic subgroups, show that they are weakly maximal,
  and that the corresponding quasi-regular representations are
  irreducible. These (infinite-dimensional) representations are
  approximated by finite-dimensional quasi-regular representations. 
  The Hecke algebras associated to these parabolic subgroups are
  commutative, so the decomposition in irreducible components of the
  finite quasi-regular representations is given by the double cosets
  of the parabolic subgroup. Since our results derive from
  considerations on finite-index subgroups, they also hold for the
  profinite completions $\widehat G$ of the groups $G$.

  The representations involved have interesting spectral properties
  investigated in~\cite{bartholdi-g:spectrum}. This paper serves as a
  group-theoretic counterpart to the studies in the mentioned paper.
  
  We study more carefully a few examples of fractal groups, and in
  doing so exhibit the first example of a torsion-free branch
  just-infinite group.

  We also produce a new example of branch just-infinite group of
  intermediate growth, and provide for it an $L$-type presentation by
  generators and relators.
\end{abstract}
\maketitle

\section{Introduction}
Fractal groups entered recently in the \emph{avant-sc\`ene} of group
theory, and are related to diverse areas such as the theory of branch
groups~\cite{grigorchuk:jibg}, automata
groups~\cite{bartholdi-g:spectrum} and so on.

Fractal groups of branch type have many interesting properties.
Namely, the first examples of groups of intermediate growth were found
in this class of groups~\cite{grigorchuk:gdegree}; the simplest
examples of infinite finitely-generated torsion groups
too~\cite{grigorchuk:burnside,gupta-s:burnside} (thus contributing
to the general Burnside problem); fractal groups provide sporadic
examples of groups of finite width with unusual associated Lie
algebra~\cite{bartholdi-g:lie}, thus answering a question by Efim
Zel'manov; etc.

It is therefore of utmost interest to pursue the study of the
algebraic, geometric and analytic properties of these groups, and in
particular their subgroup structure.

Fractal groups and branch groups are defined in the category of
profinite groups as well. These new classes of profinite groups
already started to play an important role. For instance, they gave an
answer to a question of Efim Zel'manov about groups of finite
width~\cite{bartholdi-g:lie}, they were used by Dan
Segal~\cite{segal:fgfimages} to solve in the negative a conjecture by
Alex Lubotzky, Laci Pyber and Aner Shalev~\cite{lubotzky:sg} about a
gap in the range of subgroup growths, these groups have an universal
embedding property~\cite{grigorchuk-h-z:profinite}, and it is believed
that branch groups may play an important role in Galois
theory~\cite{boston:galois}.

Fractal groups are groups acting on regular rooted trees and have
self-similarity properties inspired by those of the tree they act
on. Branch groups are groups acting on regular (or, more generally,
spherically homogeneous~\cite{grigorchuk:jibg}) rooted trees, and
having a \emph{branch structure} that endows them with properties
similar to those of the full tree automorphism group.

The action of a fractal group $G$ extends to an action on the boundary of
the tree. A \emph{parabolic} subgroup $P$ of $G$ is the stabilizer of
an element in the boundary of the tree --- or, equivalently, the
stabilizer of an infinite geodesical path starting at the root
vertex. Parabolic subgroups can be defined for any group acting on a
tree, but in the case of branch groups they have the remarkable
\emph{weak maximality} property, and the quotient spaces $G/P$
typically have polynomial growth, usually of non-integer degree.

Viewing $P$ as the stabilizer of an infinite path $e=(e_1,e_2,\dots)$,
it is approximated by the stabilizers $P_n$ of finite paths
$(e_1,\dots,e_n)$, in the sense that $P=\bigcap P_n$. The homogeneous
space $G/P$ is then also approximated by the finite spaces
$G/P_n$. These finite spaces have a limit in the Gromov sense, which
is a compact finite-dimensional space; in case its Hausdorff dimension
is not an integer, we obtain a fractal set of a new nature, as we
observed in~\cite{bartholdi-g:spectrum}. The study of such spaces is
promising.

The present paper contains several new results concerning properties
of branch fractal groups, in particular a part of their subgroup
lattice and the structure of their parabolic subgroups.

One of the main fruits of this research is the first example of
torsion-free branch just-infinite group (see
Section~\ref{sec:Gammab}).  This paper also serves as a companion
to~\cite{bartholdi-g:spectrum}, in that it studies the structure of
parabolic subgroups $P$ of fractal groups, and the decomposition of
the associated quasi-regular representations $\rho_{G/P}$.

These representations are irreducible, and that there are uncountably
many different (pairwise non-equivalent) among them. They are
infinite-dimensional, but are approximated by finite-dimensional
representations $\rho_{G/P_n}$ where $\{P_n\}$ is a sequence of
subgroups of finite index such that $P=\bigcap_{n\in\N}P_n$. For these
finite-dimensional representations we describe a decomposition in
irreducible components. This decomposition is obtained by a complete
description of the structure of the Hecke algebra associated to the
pair $(G,P)$. These Hecke algebra turn out to be abelian.

We believe branch fractal groups have good ``analytical properties''
in the sense that a sufficiently rich representation theory for these
groups, their finite images, and the corresponding profinite
completions can be developed in order to answer the main questions
about harmonic analysis on these groups --- their spectrum, the
structure of various completions of their group algebra etc.

The set $\{\rho_{G/P}|\,P\text{ is a parabolic subgroup of }G\}$ is
probably sufficiently large for this purpose, since parabolic
subgroups have the property that $\bigcap_{g\in G}P^g=1$ (for the
finite-dimensional analogue, this implies that the regular
representation $\rho_{G_n}$ is a subrepresentation of the tensor
product $\bigotimes_{|\sigma|=n}\rho_{G/\stab_G(\sigma)}$.)

We are following the first steps along this direction in the present
paper. The results given further are already used for the
computation of spectra related to fractal
groups~\cite{bartholdi-g:spectrum}, where we show that in some cases
these spectra are simple transformations of Julia sets of quadratic
maps of the complex plane.

The paper is organized as follows: in Section~\ref{sec:trees} we give
general definitions concerning groups acting on rooted trees,
introduce the congruence property, parabolic subgroups, portraits of
elements an Hausdorff dimension of closed subgroups of $\aut(\tree)$.

In Section~\ref{sec:branch} we recall the definition of branch group,
weakly branch group and regular branch group. We prove the weak
maximality of parabolic subgroups and provide a criterion evaluating
the congruence property for regular branch groups.

In Sections~\ref{sec:G}, \ref{sec:Gt}, \ref{sec:Gamma},
\ref{sec:Gammab}, \ref{sec:Gammabb} we define groups $G$, $\tilde G$,
$\Gamma$, $\overline\Gamma$, $\doverline\Gamma$, and study some
properties of these groups. We prove that $\Gamma$ and
$\overline\Gamma$ are virtually torsion-free, in contrast to $G$ and
$\doverline\Gamma$ which are torsion, and $\tilde G$ which is neither
torsion nor virtually torsion-free.

We prove that $\tilde G$ is (like $G$~\cite{grigorchuk:gdegree}
and~$\Gamma$, $\overline\Gamma$ and
$\doverline\Gamma$~\cite{bartholdi:ggs}) of intermediate growth, and
produce a presentation of $\tilde G$ which is of $L$-type, that is,
which involves finitely many relators, along with their iterates under
a word substitution. An analogous representation for $G$ was found by
Igor Lysionok~\cite{lysionok:pres}.

For each of the involved groups we draw a part of their subgroup
lattice, and provide a tree-like decomposition of their parabolic
subgroup.

We prove that $G$, $\tilde G$, $\Gamma$ and $\doverline\Gamma$ are
just-infinite branch groups, while $\overline\Gamma$ is a
just-nonsolvable weakly branch group (the first example of such a
group was given in~\cite{brunner-s-v:nonsolvable}).

Finally in Section~\ref{sec:qr} we study, for a branch group $G$, the
quasi-regular representations corresponding to parabolic subgroups $P$
and stabilizers $P_n$ of vertices at level $n$. We show that the
quasi-regular representations $\rho_{G/P}$ are irreducible, and for
the finite-dimensional representations $\rho_{G/P_n}$ we describe
their decomposition in irreducible components, which we explicit in
the case of our examples. The Hecke algebra $\Hecke(G,P_n)$, which
controls the decomposition of $\rho_{G/P_n}$ in irreducible
components, is abelian. As a consequence, the orbit structure of $P_n$
on the homogeneous $G$-space $G/P_n$ is closely related to that
decomposition.

Note that all these results --- structure of the parabolic subgroup,
lattice of finite-index subgroups, weak maximality of $P$, abelian
Hecke algebra --- hold also for the closures $\overline G$ of the
groups we consider in $\aut\tree$, which are branch fractal profinite
groups. For instance the statements about the structure of parabolic
subgroups are valid for them as well if one replaces the restricted
tree-like decomposition by an unrestricted one. Also, in our situation,
a group and its closure have the same sequences of (finite) Hecke
algebras so one can consider these algebras as associated to profinite
groups as well.  Four of our groups satisfy the congruence property,
so their closures are isomorphic to their profinite completions.

It will be important in the future to develop the theory of
representations of profinite branch groups. The results of
section~\ref{sec:qr} are a first step in this direction. As we
obtained a simple description of the double coset decomposition with
respect to a parabolic subgroup there is a hope that the classical
methods (described for instance in~\cite{curtis-r:methods}) as well as
more recent developments~\cite{michler:character} will lead to a
complete theory of the representations of the considered groups as
well as of other groups of this type.

The results in this paper are used in~\cite{bartholdi-g:spectrum}, and
are announced in~\cite{bartholdi-g:cras-parabolic}.

\subsection{Notation}
The following conventional notations shall be used: for $g,h$ in a
group $G$,
\[g^h = hgh^{-1},\qquad [g,h]=ghg^{-1}h^{-1};\]
for elements or subsets $g_1,\dots,g_n$ in $G$, the subgroup they
generate is written $\langle g_1,\dots,g_n\rangle$ and its normal
closure $\langle g_1,\dots,g_n\rangle^G$.

The symmetric subgroup on a set $\Sigma$ is written $\sym\Sigma$.

We also introduce a notation for `subsemidirect products', as follows:
\begin{defn}\label{defn:subsd}
  Let $A$ and $B$ be two subgroups of a group $G$, with $A\cap B=C$,
  and assume that $B$ is in the normalizer $N_G(A)$ and thus acts on
  $A$ by conjugation. We write $A\rtimes_C B$ for the subgroup of $G$
  generated by $A$ and $B$, and call it the \emph{subsemidirect
    product} of $A$ and $B$.
  
  If for some prime $p$ we have $C=\langle B^p,B'\rangle$, then we
  write $A\rtimes_{p-ab}B$ for the subsemidirect product of $A$ and
  $B$, and call it the \emdef{elementary abelian subsemidirect
    product} of $A$ and $B$.
\end{defn}
The motivation for this name is that there is a natural set-bijection
between $A\rtimes_CB$ and $(A\times B)/\{(c,c)\}$. The latter need
not, however, be a group, since we neither require $C$ to be normal in
$G$ nor to centralize $A$.

\section{Groups Acting on Rooted Trees}\label{sec:trees}
The groups we shall consider will all be subgroups of the group
$\aut(\tree)$ of automorphisms of a regular rooted tree $\tree$. Let
$\Sigma$ be a finite alphabet.  The vertex set of the tree
$\tree_\Sigma$ is the set of finite sequences over $\Sigma$; two
sequences are connected by an edge when one can be obtained from the
other by right-adjunction of a letter in $\Sigma$.  The top node is
the empty sequence $\emptyset$, and the children of $\sigma$ are all
the $\sigma s$, for $s\in\Sigma$.  We shall also consider the
\emdef{boundary} $\partial\tree_\Sigma$ of $\tree_\Sigma$ consisting
of the semi-infinite sequences over $\Sigma$. In most cases we shall
write $\tree$ for the rooted tree involved, when there is no ambiguity
on $\Sigma$.

We suppose $\Sigma=\Z/d\Z$, with the operation $\overline s=s+1\mod
d$.  Let $a$, called the \emdef{rooted automorphism} of
$\tree_\Sigma$, be the automorphism of $\tree_\Sigma$ defined by
$a(s\sigma)=\overline s\sigma$: it acts nontrivially on the first
symbol only, and geometrically is realized as a cyclic permutation of
the $d$ subtrees just below the root.
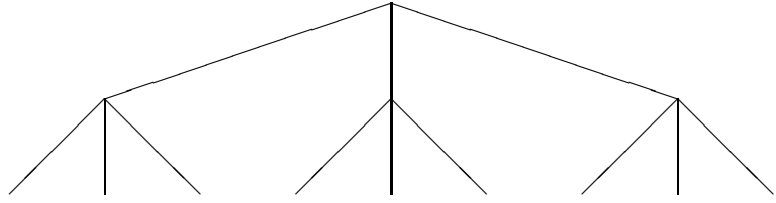
\begin{center}
\setlength{\unitlength}{3947sp}%
\begingroup\makeatletter\ifx\SetFigFont\undefined%
\gdef\SetFigFont#1#2#3#4#5{%
  \reset@font\fontsize{#1}{#2pt}%
  \fontfamily{#3}\fontseries{#4}\fontshape{#5}%
  \selectfont}%
\fi\endgroup%
\begin{picture}(4824,1224)(-11,-373)
\thinlines
{\color[rgb]{0,0,0}\put(2401,839){\line( 3,-1){1800}}
}%
{\color[rgb]{0,0,0}\put(2401,239){\line(-1,-1){600}}
}%
{\color[rgb]{0,0,0}\put(2401,239){\line( 0,-1){600}}
}%
{\color[rgb]{0,0,0}\put(4201,239){\line( 0,-1){600}}
}%
{\color[rgb]{0,0,0}\put(4201,239){\line( 1,-1){600}}
}%
{\color[rgb]{0,0,0}\put(601,239){\line( 1,-1){600}}
}%
{\color[rgb]{0,0,0}\put(601,239){\line(-1,-1){600}}
}%
{\color[rgb]{0,0,0}\put(4201,239){\line(-1,-1){600}}
}%
{\color[rgb]{0,0,0}\put(601,239){\line( 0,-1){600}}
}%
{\color[rgb]{0,0,0}\put(2401,839){\line(-3,-1){1800}}
}%
{\color[rgb]{0,0,0}\put(2401,239){\line( 1,-1){600}}
}%
{\color[rgb]{0,0,0}\put(2401,839){\line( 0,-1){600}}
}%
\end{picture}
\end{center}
Fix some $\Sigma$, and consider any subgroup $G<\aut(\tree)$. Let
$\stab_G(\sigma)$, the \emdef{vertex stabilizer} of $\sigma$, denote the
subgroup of $G$ consisting of the automorphisms that fix the sequence
$\sigma$, and let $\stab_G(n)$, the \emdef{level stabilizer}, denote the
subgroup of $G$ consisting of the automorphisms that fix all sequences
of length $n$:
\[\stab_G(\sigma)=\{g\in G|\,g\sigma=\sigma\},\qquad\stab_G(n)=\bigcap_{\sigma\in\Sigma^n}\stab_G(\sigma).\]
The $\stab_G(n)$ are normal subgroups of finite index of $G$; in
particular $\stab_G(1)$ is of index at most $d!$. Let $G_n$ be the
quotient $G/\stab_G(n)$. If $g\in\aut(\tree)$ is an automorphism fixing
the sequence $\sigma$, we denote by $g_{|\sigma}$ the element of
$\aut(\tree)$ corresponding to the restriction to sequences starting
by $\sigma$:
\[\sigma g_{|\sigma}(\tau)=g(\sigma\tau).\]
As the subtree starting from any vertex is isomorphic to the initial
tree $\tree_\Sigma$, we obtain this way a map
\begin{equation}\label{eq:phi}
  \phi:\begin{cases}\stab_{\aut(\tree)}(1)\to\aut(\tree)^\Sigma\\
    h\mapsto (h_{|0},\dots,h_{|d-1})\end{cases}
\end{equation}
which is an embedding.  For a sequence $\sigma$ and an automorphism
$g\in\aut(\tree)$, we denote by $g^\sigma$ the element of
$\aut(\tree)$ acting as $g$ on the sequences starting by $\sigma$, and
trivially on the others:
\[g^\sigma(\sigma\tau)=\sigma g(\tau),\qquad g^\sigma(\tau)=\tau\text{ if }\tau\text{ doesn't start by }\sigma.\]
The \emdef{rigid stabilizer} of $\sigma$ is
$\rist_G(\sigma)=\{g^\sigma|\,g\in G\}\cap G$. We also set
\[\rist_G(n)=\langle\rist_G(\sigma)|\,\sigma\in\Sigma^n\rangle=\prod_{|\sigma|=n}\rist_G(\sigma)\]
and call it the \emdef{rigid level stabilizer} ($\prod$ denotes direct
product).  We say $G$ has \emdef{infinite rigid stabilizers} if all
the $\rist_G(\sigma)$ are infinite.

\begin{defn}
  A subgroup $G<\aut(\tree)$ is \emph{spherically transitive} if the
  action of $G$ on $\Sigma^n$ is transitive for all $n\in\N$.

  $G$ is \emdef{fractal} if for every vertex $\sigma$ of
  $\tree_\Sigma$ one has $\stab_G(\sigma)_{|\sigma}\cong G$, where the
  isomorphism is given by identification of $\tree_\Sigma$ with its
  subtree rooted at $\sigma$.

  $G$ has the \emdef{congruence property} if every finite-index
  subgroup of $G$ contains $\stab_G(n)$ for some $n$ large enough.
\end{defn}

\begin{lem}\label{lem:fractal}
  The group $G<\aut(\tree)$ is fractal if and only if
  $\phi_{|\stab_G(1)}:\stab_G(1)\to\aut(\tree)^\Sigma$ is a subdirect
  embedding into $G\times\dots\times G$, i.e.\ if it is an embedding
  that is surjective on each factor.
\end{lem}
\begin{proof}
  If $p_i\phi_{|G}\neq G$ for some projection $p_i$ on the vertex $i$, 
  then $\stab_G(i)_{|i}\neq G$ so $G$ is not fractal. We now suppose
  $\phi_{|G}$ is a subdirect embedding and prove by induction that
  $\stab_G(\sigma)_{|\sigma}\cong G$ for all $\sigma$.

  The induction basis, for $|\sigma|=1$, is equivalent to the
  hypothesis. Now by induction $G\to G^{\Sigma^{n-1}}$ is a subdirect
  embedding, and each factor $G$ maps to $G^\Sigma$ by
  $\phi_{|G}$. The composition of two subdirect embeddings is again
  subdirect, so $G\to G^{\Sigma^n}$ is subdirect.
\end{proof}
For fractal groups, we usually shall write $\phi$ instead of
$\phi_{|G}$.

\begin{lem}
  A fractal group is spherically transitive if and only if it acts
  transitively on the first level $\Sigma$.
\end{lem}
\begin{proof}
  Assume by induction that $G$ is fractal and transitive on
  $\Sigma^{n-1}$, the induction starting at $n=2$. Since $\phi$ is
  subdirect, $G$ is transitive on $i\Sigma^{n-1}$ for all
  $i\in\Sigma$, and since it is transitive on $\Sigma$, it is also
  transitive on $\Sigma^n$.
\end{proof}

The full automorphism group $\aut(\tree)$ has the structure of a
profinite group: it is approximated by the finite groups
$\aut(\tree)_n=\aut(\tree)/\stab(n)$, and we have
\[\aut(\tree)=\varprojlim_{n\to\infty}\aut(\tree)_n.\]
More on the topology of $\aut(\tree)$ is said
in~\cite{bartholdi-g:spectrum}.  The following lemma follows directly
from the definition of a profinite completion:
\begin{lem}
  Let $G<\aut(\tree)$ have the congruence property. Then its profinite
  completion $\widehat G$ is isomorphic (as a profinite group) to its
  closure $\overline G$ in $\aut(\tree)$. If moreover $G$ is contained
  in a Sylow pro-$p$-subgroup $\aut_p(\tree)$, then $\overline G$
  is isomorphic to the pro-$p$ completion $\widehat G_p$ of $G$.
\end{lem}
\begin{proof}
  By the congruence property, $\{\overline{\stab_G(n)}\}$ is a basis of
  neighbourhoods of the identity in $\widehat G$.
\end{proof}

\begin{defn}
  Assume $G<\aut(\tree)$ is given, with a subset $S\subset G$.
  The \emdef{portrait} of $g\in G$ with respect to $S$ is a subtree of
  $\tree$, with inner vertices labeled by $\sym\Sigma$ and leaf
  vertices labeled by $S\cup\{1\}$. It is defined recursively as follows: if
  $g\in S\cup\{1\}$, the portrait of $g$ is the subtree reduced to the root
  vertex, labeled by $g$ itself. Otherwise, let $\alpha\in\sym\Sigma$
  be the permutation of the top branches of $\tree$ such that
  $g\alpha^{-1}\in\stab_G(1)$; let
  $(g_0,\dots,g_{d-1})=\phi(g\alpha^{-1})$ and let $\tree_i$ be the
  portrait of $g_i$. Then the portrait of $g$ is the subtree of
  $\tree$ with $\alpha$ labeling the root vertex and subtrees
  $\tree_0,\dots,\tree_{d-1}$ connected to the root.

  The \emdef{portrait} of $g\in G$ is its portrait with respect to
  $\emptyset$. The element $g$ is called \emdef{finitary} if its
  portrait is finite.

  The \emdef{depth} of $g\in G$ is the height (length of a maximal
  path starting at the root vertex) $\partial(g)\in\N\cup\{\infty\}$
  of the portrait of $g$.
\end{defn}

We now suppose $d=p$ is prime, and consider $\aut_*(\tree)$, the Sylow
pro-$p$ subgroup of $\aut(\tree)$ consisting of all elements $g$ whose
portrait is labeled by powers of the cycle $(0,1,\dots,d-1)$. It has
the structure of an infinitely iterated wreath product
\[\aut_*(\tree)=\Z/p\Z\wr\Z/p\Z\wr\dots.\]
For a closed subgroup $G$ of $\aut_*(\tree)$, its \emdef{Hausdorff
  dimension} $\dim_*(G)$ is defined in~\cite{barnea-s:hausdorff} as
\[\dim_*(G)=\liminf_{n\to\infty}\frac{p-1}{p^n}\log_p|G_n|=\liminf_{n\to\infty}\frac{\log_p|G_n|}{\log_p|\aut_*(\tree)_n|}.\]
In particular, the Hausdorff dimension of $\aut_*(\tree)$ is $1$, and
$\dim_*$ is invariant upon taking finite-index subgroups.

\section{Branch Groups}\label{sec:branch}
We consider now a special class of groups acting on rooted trees. We
shall always implicitly assume they act spherically transitively.
\begin{defn}\label{defn:b}\begin{enumerate}
  \item $G$ is a \emdef{regular branch} group if it has a finite-index
    normal subgroup $K<\stab_G(1)$ such that
    \[K^\Sigma<\phi(K).\]\label{defn::rb}
    It is then said to be \emdef{regular branch over $K$}.
  \item A subgroup $G<\aut(\tree)$ is a \emdef{branch group} if for
    every $n\ge1$ the subgroup $\rist_G(n)$ has finite index in
    $G$.\label{defn::b}
  \item $G$ is a \emdef{weakly branch} group if all of its rigid
    stabilizers $\rist_G(\sigma)$ are non-trivial.\label{defn::wb}
  \end{enumerate}
\end{defn}
Note that the definition of a branch group admits
an even more general setting --- see~\cite{grigorchuk:jibg}. Four of
our examples will be regular branch groups, and the last one will not
be a branch, but rather a weakly branch group. The following lemma shows
that, for fractal groups, \ref{defn::rb} implies \ref{defn::b} implies
\ref{defn::wb} in Definition~\ref{defn:b}.
\begin{lem}\label{lem:rb->b}
  If $G$ is a fractal, regular branch group, then it is a branch group.
  If $G$ is a branch group, then it is a weakly branch group.
\end{lem}
\begin{proof}
  Assume $G$ is a regular branch group over its subgroup $K$. Clearly
  $K^{\Sigma^n}$ can be viewed, through $\phi^n$, as a subgroup of
  $\rist_G(n)$, and is of finite index in $G^{\Sigma^n}$, so
  $\rist_G(n)$ is of finite index in $G$. The second implication holds
  because branch groups are infinite, and `finite index in an infinite
  group' is stronger than `non-trivial'.
\end{proof}

Note also that if all rigid stabilizers are non-trivial, then they are
all infinite; moreover,
\begin{lem}\label{lem:inforbits}
  Let $G$ be a weakly branch group, $\sigma\in\Sigma^n$ a vertex, and
  $\overline\sigma=\sigma_1\dots\sigma_n\sigma_{n+1}\dots\in\Sigma^\N$
  and infinite ray extending $\sigma$. Then the
  $\rist_G(\sigma)$-orbit of $\overline\sigma$ is infinite.
\end{lem}
\begin{proof}
  It suffices to show that the orbit of
  $v_k=\sigma_1\dots\sigma_{n+k}$ becomes arbitrarily large as $k$
  increases.  Since $\rist_G(\sigma)$ is non-trivial, it contains $g$
  moving $\sigma\tau$ to $\sigma\tau'$ for some
  $\tau,\tau'\in\Sigma^k$.  Since $G$ is spherically transitive, it
  contains $h$ moving $v_k$ to $\sigma\tau$. Consider now $g^h$: it
  belongs to $\rist_G(\sigma)$, and does not fix $v_k$, whence $v_k$'s
  orbit contains at least $2$ points.
  
  The argument applied to $v_k$ shows that some $v_{k+k'}$ has least
  $2$ points in its $\stab_G(v_k)$-orbit, so at least $4$ points in
  its $\stab_G(\sigma)$-orbit; and this process can be repeated an
  arbitrary number of times to produce vertices $v_{k+\dots+k^{(j)}}$
  with at least $2^{j+1}$ points in their orbit.
\end{proof}

If $G$ is a regular branch group over its subgroup $K$, the following
notations is also introduced: given a subgroup $L$ of $K$, we write
$L_{(0)}=L$ and inductively $L_{(n)}=\phi^{-1}(L_{(n-1)}^\Sigma)$.
These $L_{(i)}$ form a sequence of subgroups of $L$ with
$\bigcap_{n\ge0}L_{(n)}=\{1\}$.

\begin{defn}
  A group $G$ is \emdef{just-infinite} if it is infinite, and for any
  non-trivial normal subgroup $N$ the quotient $G/N$ is finite.
\end{defn}
Note that in checking just-infiniteness one may restrict one's
attention to subgroups $N=\langle g\rangle^G$, i.e.\ normal closures
of a non-trivial element of $G$. We will use the following criterion:
\begin{prop}\label{prop:jinf}
  Let $G$ be regular branch over $K$. Then $G$ is just infinite if
  and only if $|K:K'|<\infty$.
\end{prop}
\begin{proof}
  Clearly if $K'$ is of infinite index in $K$ then $\langle
  K'\rangle^G$ is of infinite index in $G$, and is not trivial ($K$
  clearly cannot be abelian) so $G$ is not just infinite.
  
  Conversely, assume $|K:K'|<\infty$ and let $G\ni g\neq1$. Let
  $N=\langle g\rangle^G$; we will show that $N$ is of finite index.
  Determine $n$ such that
  $g\in\stab_{\aut(\tree)}(n-1)\setminus\stab_{\aut(\tree)}(n)$. Then
  there is a sequence $\sigma$ of length $n-1$ such that
  $g_{|\sigma}\not\in\stab_{\aut(\tree)}(1)$.  Choose now two elements
  $k_1,k_2$ of $K$.  Because $G$ is branched on $K$, it contains for
  $i=1,2$ the elements $\xi_i=k_i^{\sigma0}$.  Let $\eta=[\xi_1,g]\in
  N$. It fixes all sequences except: those starting by $\sigma0$, upon
  which it acts as $k_1$, and possibly those starting by $\sigma x$
  for $x\ge1$. Consider $\zeta=[\eta,\xi_2]\in N$. Clearly
  $\zeta=[k_0,k_1]^{\sigma0}$; as the commutator $[k_0,k_1]$ was
  chosen arbitrarily, it follows that $N$ contains ${K'}^{\sigma0}$;
  and as $N$ is normal, it contains $K'\times\dots\times K'$, the
  product having $d^n$ factors. Now $K'$ is of finite index in $K$
  which is of finite index in $G$, so $K'\times\dots\times K'$ is of
  finite index in $G$ and the same holds for $N$.
\end{proof}

\begin{defn}
  A group $G$ is \emdef{just-non-solvable} if it is not solvable, but
  all its non-trivial quotients are solvable.
\end{defn}
\begin{prop}\label{prop:jns}
  Let $G$ be regular branch over $K$. Then $G$ is just-non-solvable if
  and only if $G/K_{(1)}$ is solvable. In particular, if $d$ is prime,
  every regular branch subgroup of $\aut_*(\tree)$ is
  just-non-solvable.
\end{prop}
\begin{proof}
  Let $G$ be just-non-solvable. Then $K_{(1)}$ is a non-trivial normal
  subgroup, so $G/K_{(1)}$ is solvable.

  Let now $N$ be a non-trivial normal subgroup of $G$. It is shown
  in~\cite[Theorem~4]{grigorchuk:jibg} that $N$ contains $K'_{(n)}$ for
  some $n\in\N$, so it suffices to show that all $G/K'_{(n)}$ are
  solvable. Consider the chain
  \[G/K'_{(n)}\triangleright K/K'_{(n)}\triangleright\dots\triangleright
  K_{(n)}/K'_{(n)}\cong K/K'\times\dots\times K/K'.\] The last group is
  abelian, hence solvable ; successive quotients in the sequence are
  also solvable because
  \[(K_{(i)}/K'_{(n)})\big/(K_{(i+1)}/K'_{(n)}) = K_{(i)}/K_{(i+1)} \cong
  K/K_{(1)}\times\dots\times K/K_{(1)},\] and $K/K_{(1)}$ is solvable
  by assumption. Also, $(G/K'_{(n)})/(K/K'_{(n)})\cong G/K$ is
  solvable; therefore $G/K'_{(n)}$ is an extension of solvable groups,
  so is solvable.
  
  In case $d$ is prime and $G$ is a branch subgroup of
  $\aut_*(\tree)$, the quotient $G/K_{(1)}$ is a finite $d$-group, so
  is solvable.
\end{proof}

The following criterion describes which branch groups enjoy the
congruence property:
\begin{prop}
  Let $G$ be regular branch over $K$. Then $G$ has the congruence
  property if $K'$ contains $\stab_G(m)$ for some $m\in\N$.
\end{prop}
\begin{proof}
  Let $N$ be a finite-index subgroup of $G$. By replacing $N$ with its
  core $\bigcap_{g\in G}N^g$, still of finite index, we may suppose
  $N$ is normal in $G$. By~\cite[Theorem~4]{grigorchuk:jibg}, $N$
  contains $K'_{(n)}$ for some $n\in\N$, so
  \[N > K'_{(n)} > \stab_G(m)_{(n)} > \stab_G(m+n).\]
\end{proof}

As an example of regular branch group not enjoying the congruence
property, consider $G=\aut_f(\tree)$, the automorphisms of $\tree$
whose action $\phi_v\in\sym\Sigma$ is non-trivial only at finitely
many vertices $v$, and its subgroup
$H = \{g\in G|\,\prod_{v\in\Sigma^*}\phi_v\in\mathfrak{Alt}_\Sigma\}.$
Here $G$ is regular branch, with $K=G$, but $H$ does not contain any
level stabilizer.

When furthermore $G$ is finitely generated, the following
`quantitative congruence property' shall be useful to prove equalities
among subgroups:
\begin{prop}[Quantitative Congruence Property]\label{prop:qcp}
  Let $G$ be regular branch over $K$, finitely generated by the set
  $S$ and with the congruence property. Let $n$ be minimal such that
  $\langle s\rangle^G\ge\stab_G(n)$ for all generators $s\in S$. Let $m$
  be minimal such that $K'$ contains $\stab_G(m)$.
  
  Let $N$ be any non-trivial normal subgroup of $G$ and $1\neq g\in
  N$. Then $N$ contains $\stab_G(\partial(g)+m+n)$.
\end{prop}
\begin{proof}
  This follows again from Theorem~4 in~\cite{grigorchuk:jibg}.
\end{proof}

\subsection{Parabolic subgroups} In the context of groups acting on
the hyperbolic space, a parabolic subgroup is the stabilizer of a
point on the boundary. We give here a few general facts concerning
parabolic subgroups of fractal or branched groups, and recall some
results on growth of groups and sets on which they act.

\begin{defn}
  Let $\tree=\Sigma^*$ be a rooted tree. A \emdef{ray} $e$ in $\tree$
  is an infinite geodesic starting at the root of $\tree$, or
  equivalently an element of $\partial\tree=\Sigma^\N$.
  
  Let $G<\aut(\tree)$ be any subgroup and $e$ be a ray. The associated
  \emdef{parabolic subgroup} is $P_e=\stab_G(e)$.
\end{defn}

The following important facts are easy to prove:
\begin{itemize}
\item $\bigcap_{e\in\partial\tree}P_e=\bigcap_{g\in G}P^g=1$.
\item Let $e$ be an infinite ray and define the subgroups
  $P_n=\stab_G(e_1\dots e_n)$. Then the $P_n$ have index $d^n$ in $G$
  (since $G$ acts transitively) and satisfy
  \[P_e=\bigcap_{n\in\N}P_n.\]
\item $P$ has infinite index in $G$, and has the same image as $P_n$
  in the quotient $G_n=G/\stab_G(n)$.
\end{itemize}

\begin{defn}\label{defn:growth}
  Let $G$ be a group generated by a finite set $S$, let $X$ be a set
  upon which $G$ acts transitively, and choose $x\in X$. The
  \emdef{growth} of $X$ is the function $\gamma:\N\to\N$ defined by
  $$\gamma(n)=|\{gx\in X|\,|g|\le n\}|,$$
  where $|g|$ denotes the minimal length of $g$ when written as a word
  over $S$.  By the \emdef{growth of $G$} we mean the growth of the
  action of $G$ on itself by left-multiplication.
  
  Given two functions $f,g:\N\to\N$, we write $f\preceq g$ if there is
  a constant $C\in\N$ such that $f(n)<Cg(Cn+C)+C$ for all $n\in\N$,
  and $f\sim g$ if $f\preceq g$ and $g\preceq f$. The equivalence
  class of the growth of $X$ is independent of the choice of $S$ and
  of $x$.

  $X$ is of \emdef{polynomial growth} if $\gamma(n)\preceq n^d$ for some
  $d$. It is of \emdef{exponential growth} if $\gamma(n)\succeq
  e^n$. It is of \emdef{intermediate growth} in the remaining
  cases. This trichotomy does not depend on the choice of $x$ or $S$.
\end{defn}

\begin{defn}
  Two infinite sequences $\sigma,\tau:\N\to\Sigma$ are
  \emdef{confinal} if there is an $N\in\N$ such that $\sigma_n=\tau_n$
  for all $n\ge N$.

  Confinality is an equivalence relation, and equivalence classes are
  called \emdef{confinality classes}.
\end{defn}

The following result is due to Volodymyr Nekrashevych and Vitaly
Sushchansky.
\begin{prop}\label{prop:confinal}
  Let $G$ be a group acting on a regular rooted tree $\tree$, and
  assume that for any generator $g\in G$ and infinite sequence
  $\sigma$, the sequences $\sigma$ and $g\sigma$ differ only in
  finitely many places. Then the confinality classes of the action of
  $G$ on $\partial\tree$ are unions of orbits. If moreover
  $\stab_G(\sigma)$ contains the rooted automorphism $a$ for all
  $\sigma\in\tree$, the orbits of the action are confinality classes.
\end{prop}

The proof of the following result appears in~\cite{bartholdi-g:spectrum}.
\begin{prop}\label{prop:polygrowth}
  Let $G<\aut(\tree)$ satisfy the conditions of
  Proposition~\ref{prop:confinal}, and suppose that for the map
  $\phi:g\mapsto(g_1,\dots,g_d)$ defined in~(\ref{eq:phi}) there are
  constants $\lambda,\mu$ such that $|g_i|\le\lambda|g|+\mu$ for all
  $i$. Let $P$ be a parabolic subgroup. Then $G/P$, as a $G$-set, is
  of polynomial growth of degree at most $\log_{1/\lambda}(d)$.  If
  moreover $G$ is spherically transitive, then $G/P$'s asymptotical
  growth is polynomial of degree $\log_{1/\lambda'}(d)$, with
  $\lambda'$ the infimum of the $\lambda$ as above.
\end{prop}

\begin{defn}
  Let $G$ be a branch group, and $H$ any subgroup. We say $H$ is
  \emdef{weakly maximal} if $H$ is of infinite index in $G$, but
  all subgroups of $G$ strictly containing $H$ are of finite index in $G$.
\end{defn}
Note that every infinite finitely generated group admits maximal
subgroups, by Zorn's lemma.

Note also that some branch groups may not contain any infinite-index
maximal subgroup; this is the case for $G$, as was shown by Ekaterina
Pervova.

\begin{prop}
  Let $P$ be a parabolic subgroup of a regular branch group
  $G$. Then $P$ is weakly maximal.
\end{prop}
\begin{proof}
  Assume $G$ is regular branch over $K$, and $P=\stab_G(e)$. Recall that
  $G$ contains a product of $d^n$ copies of $K$ at level $n$, and
  clearly $P$ contains a product of $d^n-1$ copies of $K$ at level
  $n$, namely all but the one indexed by the vertex $e_1\dots e_n$.
  
  Take $g\in G\setminus P$. There is then an $n\in\N$ such that
  $g(e_n)\neq e_n$, so $\langle P,P^g\rangle$ contains the product of
  $d^n$ copies of $K$ at level $n$, hence is of finite index in $G$.
\end{proof}

\section{The Group $G$}\label{sec:G}
We give here the basic facts we will use about the first of
Grigorchuk's examples, the group
$G$~\cite{grigorchuk:burnside,harpe:ggt}. We apply the discussion
of the previous section to $\Sigma=\{0,1\}$.  Recall $a$ is the
automorphism permuting the top two branches of $\tree_2$. Define
recursively by $b$ the automorphism acting as $a$ on the right branch
and $c$ on the left, by $c$ the automorphism acting as $a$ on the
right branch and $d$ on the left, and by $d$ the automorphism acting
as $1$ on the right branch and $b$ on the left. In formul\ae,
\begin{alignat*}{2}
  b(0x\sigma)&=0\overline x\sigma,&\qquad b(1\sigma)&=1c(\sigma),\\
  c(0x\sigma)&=0\overline x\sigma,&\qquad c(1\sigma)&=1d(\sigma),\\
  d(0x\sigma)&=0x\sigma,&\qquad d(1\sigma)&=1b(\sigma).
\end{alignat*}
$G$ is the group generated by $\{a,b,c,d\}$.  It is readily checked
that these generators are of order $2$ and that $\{1,b,c,d\}$
constitute a Klein group; one of the generators $\{b,c,d\}$ can thus be
omitted.

We write $H_n=\stab_G(n)$ and $H=H_1$. Explicitly, the map $\phi$
restricts to
\[\phi:\begin{cases}H\to G\times G\\
  b\mapsto(a,c),\quad b^a\mapsto(c,a)\\ c\mapsto(a,d),\quad
  c^a\mapsto(d,a)\\ d\mapsto(1,b),\quad d^a\mapsto(b,1).\end{cases}\]
Define also the following subgroups of $G$:
\begin{align*}
  B&=\langle b\rangle^G = \langle b, b^a, b^{da}, b^{ada}\rangle,\\
  K&=\langle(ab)^2\rangle^G = \langle (ab)^2, (abad)^2, (adab)^2\rangle.
\end{align*}

The group $G$
\begin{itemize}
\item is an infinite torsion $2$-group.
\item is of intermediate growth.
\item is amenable.
\item is fractal and branched on its subgroup $K$.
\item is just-infinite.
\item is residually finite.
\item has an infinite recursive presentation~\cite{lysionok:pres} of
  $L$-type
  \[G = \langle a,b,c,d|\,a^2,b^2,c^2,d^2,bcd,\sigma^i(ad)^4,\sigma^i(adacac)^4\,(i\in\N)\rangle,\]
  where $\sigma$ is the substitution on $\{a,b,c,d\}^*$ defined by
  \[\sigma(a)=aca,\quad\sigma(b)=d,\sigma(c)=b,\sigma(d)=c,\]
  which induces an injective expanding endomorphism of $G$ of
  infinite-index image.  Moreover none of the relators of $G$ are
  superfluous~\cite{grigorchuk:bath}.
\item The subgroup $B$ is of index $8$ in $G$ and $K$ is of index
  $16$.  Also, $K$ contains $\stab_G(3)$ and $K'$ contains $K_{(2)}$, so $G$
  has the congruence property.
\item The quotients $G_n=G/\stab_G(n)$ have order $2^{5\cdot2^{n-3}+2}$
  for $n\ge3$ (and order $2^{2^n-1}$ for $n\le3$). Therefore the
  closure $\overline G$ of $G$ in $\aut(\tree)$ is isomorphic to the
  profinite completion $\widehat G$ and is a pro-$2$-group. It has
  Hausdorff dimension $5/8$~\cite{grigorchuk:jibg}.
\end{itemize}

Most of these facts are known, and appear in the extensive
reference~\cite{harpe:ggt} and in~\cite{grigorchuk:jibg}.  One can
then check by direct computation that $K$ is of index $16$.  To prove
that $G$ is regular branch, set $x=(ab)^2$. Then one has
$\phi([x,d])=(x,1)$ so by conjugation $\phi(K)>K\times1$ and thus
$\phi(K)>K\times K$.  By direct computation, $K'$ is of index $64$ in
$K$, so $G$ is just-infinite.

For all other computations, we propose an alternate method of proof,
based on the following
\begin{lem}\label{lem:Gqcp}
  $G$ satisfies the Quantitative Congruence Property for $m=4$ and
  $n=3$.
\end{lem}
\begin{proof}
  This follows from the above description of $K$.
\end{proof}

\begin{prop} We have
  \begin{align*}
    \phi(H)&=(B\times B)\rtimes\langle\phi(c),\phi(c^a)\rangle,\\
    \phi(B)&=(K\times K)\rtimes_{\langle\phi(ab)^8\rangle}\langle\phi(b),\phi(b^a)\rangle,\\
    \phi(K)&=(K\times K)\rtimes_{\langle\phi(ab)^8\rangle}\langle\phi(ab)^2\rangle,
  \end{align*}
  with the notation introduced in Definition~\ref{defn:subsd}.
\end{prop}
\begin{proof}
  Each of these subgroups $H,B,K$ are normal finite-index subgroups of
  $G$. By the Quantitative Congruence Property, they are all contained
  in some $\stab_G(n)$. It is therefore equivalent to study them
  directly or to study their images in $G_n=G/\stab_G(n)$, which is a
  finite group. A computer algebra system, like~\cite{gap:manual}, can
  then be used to derive their structure.
\end{proof}

\subsection{Low-Index subgroups}
$G$ has $7$ subgroups of index $2$:
\begin{alignat*}{3}
  \langle b,ac\rangle,&&\langle c,ad\rangle,&&\langle d,ab\rangle,\\
  \langle b,a,a^c\rangle,&&\langle c,a,a^d\rangle,&&\langle d,a,a^b\rangle,\\
  &&H = \langle b,c,b^a,c^a\rangle.
\end{alignat*}
As can be computed from its presentation~\cite{lysionok:pres} and a
computer algebra system~\cite{gap:manual}, $G$ has the following
subgroup count:

\centerline{\begin{tabular}{c|cc|cc|}
  Index & Subgroups & Normal & In $H$ & Normal\\
  \hline
  1 & 1 & 1 & 0 & 0\\
  2 & 7 & 7 & 1 & 1\\
  4 & 19 & 7 & 9 & 4\\
  8 & 61 & 7 & 41 & 7\\
  16 & 237 & 5 & 169 & 5\\
  32 & 843 & 3 & 609 & 3\\
  \hline
\end{tabular}}

\subsection{Normal closures of generators} They are as follows:
\begin{align*}
  A = \langle a\rangle^G &= \langle a,a^b,a^c,a^d\rangle,&G/A&\cong\Z/2\Z\times\Z/2\Z,\\
  B = \langle b\rangle^G &= \langle b,b^a,b^{ad},b^{ada}\rangle,& G/B &\cong D_8,\\
  C = \langle c\rangle^G &= \langle c,c^a,c^{ad},c^{ada}\rangle,& G/C &\cong D_8,\\
  D = \langle d\rangle^G &= \langle d,d^a,d^{ac},d^{aca}\rangle,& G/D &\cong D_{16}.
\end{align*}

\subsection{Some other subgroups} To complete the picture, we
introduce the following subgroups of $G$:
\begin{align*}
  K &= \langle (ab)^2\rangle^G,&
  L &= \langle (ac)^2\rangle^G,&
  M &= \langle (ad)^2\rangle^G,\\
  \overline B &= \langle B,L\rangle,&
  \overline C &= \langle C,K\rangle,&
  \overline D &= \langle D,K\rangle,\\
  && T = K^2 &= \langle (ab)^4\rangle^G,\\
  T_{(m)} &= \underbrace{T\times\dots\times T}_{2^m},&
  K_{(m)} &= \underbrace{K\times\dots\times K}_{2^m},&
  N_{(m)} &= T_{(m-1)}K_{(m)}.
\end{align*}

\begin{figure}
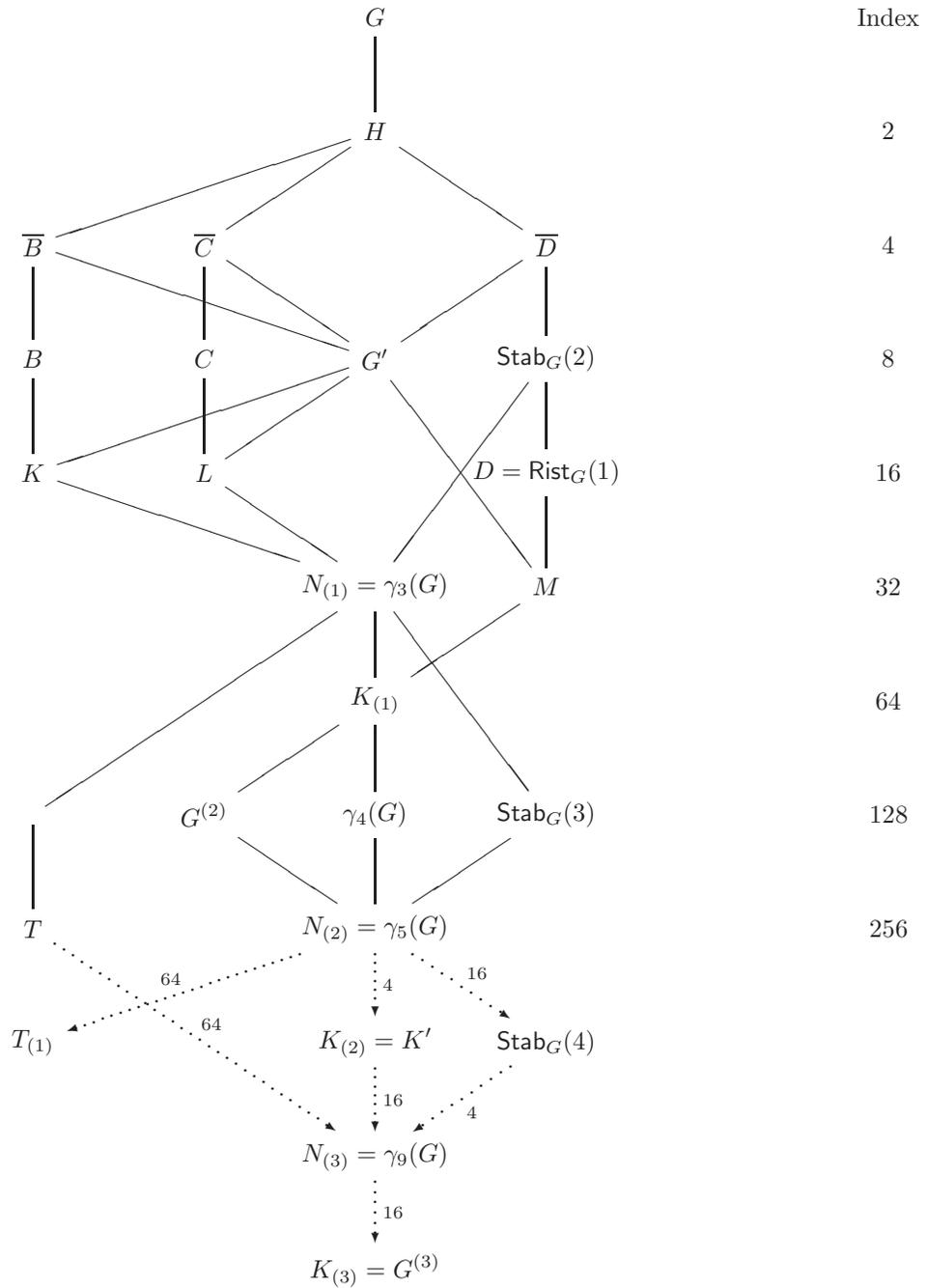

  \[\begin{diagram}\dgARROWLENGTH=1.8em
    \node[3]{G} \arrow{s,-} \node[3]{\text{Index}}\\
    \node[3]{H} \arrow{wsw,-} \arrow{sw,-} \arrow{se,-} \node[3]{2}\\
    \node{\overline B} \arrow{s,-} \arrow{ese,-}
    \node{\overline C} \arrow{s,-} \arrow{se,-}
    \node[2]{\overline D} \arrow{s,-} \arrow{sw,-} \node[2]{4}\\
    \node{B} \arrow{s,-}
    \node{C} \arrow{s,-}
    \node{G'} \arrow{wsw,-} \arrow{sw,-} \arrow{sse,-}
    \node{\stab_G(2)} \arrow{s,-} \arrow{ssw,-} \node[2]{8}\\
    \node{K} \arrow{ese,-}
    \node{L} \arrow{se,-}
    \node[2]{D=\rist_G(1)} \arrow{s,-} \node[2]{16}\\
    \node[3]{N_{(1)}=\gamma_3(G)} \arrow{s,-} \arrow{sse,-} \arrow[2]{sw,-}
    \node{M} \arrow{sw,-} \node[2]{32}\\
    \node[3]{K_{(1)}} \arrow{sw,-} \arrow{s,-} \node[3]{64}\\
    \node{} \arrow{s,-}
    \node{G^{(2)}} \arrow{se,-}
    \node{\gamma_4(G)} \arrow{s,-}
    \node{\stab_G(3)} \arrow{sw,-} \node[2]{128}\\
    \node{T} \arrow[2]{se,l,..}{64}
    \node[2]{N_{(2)}=\gamma_5(G)} \arrow{wsw,l,..}{64}
    \arrow{s,r,..}{4} \arrow{se,t,..}{16} \node[3]{256}\\
    \node{T_{(1)}}
    \node[2]{K_{(2)}=K'} \arrow{s,r,..}{16}
    \node{\stab_G(4)} \arrow{sw,r,..}{4}\\
    \node[3]{N_{(3)}=\gamma_9(G)} \arrow{s,r,..}{16}\\
    \node[3]{K_{(3)}=G^{(3)}}
  \end{diagram}\]
  \caption{The top of the lattice of normal subgroups of $G$ below
    $H$. The index of the inclusions are indicated next to the edges.}
  \label{fig:Glattice}
\end{figure}

\begin{thm}
  \begin{itemize}
  \item In the Lower Central Series, $\gamma_{2^m+1}(G) = N_{(m)}$ for all
    $m\ge 1$.
  \item In the Derived Series, $K^{(n)}=\rist_G(2n)$ for all $n\ge2$
    and $G^{(n)}=\rist_G(2n-3)$ for all $n\ge3$.
  \item The rigid stabilizers satisfy
    \[\rist_G(n)=\begin{cases} D & \text{ if }n=1,\\
      K_{(n)} & \text{ if }n\ge2.
    \end{cases}\]
  \item The level stabilizers satisfy
    \[\stab_G(n)=\begin{cases} H & \text{ if }n=1,\\
      \langle D,T\rangle & \text{ if }n=2,\\
      \langle N_{(2)},(ab)^4(adabac)^2\rangle & \text{ if }n=3,\\
      \underbrace{\stab_G(3)\times\dots\times\stab_G(3)}_{2^{n-3}} & \text{ if }n\ge4.
    \end{cases}\]
  \item There is for all $\sigma\in\Sigma^n$ a surjection
    $\cdot_{|\sigma}:\stab_G(n)\twoheadrightarrow G$ given by projection
    on the factor indexed by $\sigma$.
  \end{itemize}
  The top of the lattice of normal subgroups of $G$ below $H$ is given
  in Figure~\ref{fig:Glattice}.
\end{thm}
\begin{proof}
  The first three points are proven by Alexander Rozhkov
  in~\cite{rozhkov:habilitation}. To prove the fourth assertion, we
  apply Lemma~\ref{lem:Gqcp} to determine the structure of $\stab_G(n)$
  for $n\le4$, and note that $\stab_G(4)=\stab_G(3)\times\stab_G(3)$.
  
  For the last statement, note that $s=(ab)^4(adabac)^2$ belongs to
  $\stab_G(3))$, and that
  $\phi^n(\sigma^{n-3}(s))=(1,\dots,1,ba,d,d,ba,a,c,a,c)$ giving,
  after conjugation, any generator of $G$ in any position on any level
  $n$.
\end{proof}

\subsection{The Subgroup $P$}\label{subs:subP}
Let $e$ be the ray $1^\infty$ and let $P$ be the corresponding
parabolic subgroup.

\begin{thm}\label{thm:decomP}
  $P/P'$ is an infinite elementary $2$-group generated by the images of
  $c$, $d=(1,b)$ and of all elements of the form $(1,\dots,1,(ac)^4)$
  in $\rist_G(n)$ for $n\in\N$. The
  following decomposition holds:
  \[P = \bigg(B\times \Big(\big(K\times ((K\times\dots)\rtimes\langle (ac)^4\rangle)\big)\rtimes\langle b,(ac)^4\rangle\Big)\bigg)\rtimes\langle c,(ac)^4\rangle,\]
  where each factor (of nesting $n$) in the decomposition acts on the
  subtree just below some $e_n$ but not containing $e_{n+1}$.
\end{thm}
Note that we use the same notation for a subgroup $B$ or $K$ acting on
a subtree, keeping in mind the identification of a subtree with the
original tree. The same convention will hold for
Theorems~\ref{thm:decomPtilde}, \ref{thm:decomPG}, \ref{thm:decomPGB},
\ref{thm:decomPGBB}, and all related propositions. Note also that
$\phi$ is omitted when it would make the notations too heavy.

\begin{proof}
  Define the following subgroups of $G_n$:
  \begin{gather*}
    B_n = \langle b\rangle^{G_n};\qquad K_{(n)} = \langle (ab)^2\rangle^{G_n};\\
    Q_n = B_n\cap P_n;\qquad R_n = K_{(n)}\cap P_n.
  \end{gather*}
  Then the theorem follows from the following proposition.
\end{proof}

\begin{figure}
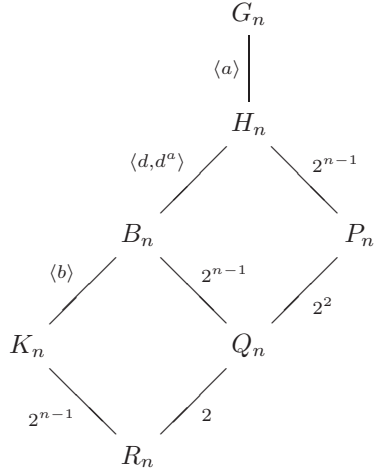

  \[\begin{diagram}
    \node[3]{G_n}\arrow{s,l,-}{\langle a\rangle}\\
    \node[3]{H_n}\arrow{sw,l,-}{\langle d,d^a\rangle}\arrow{se,l,-}{2^{n-1}}\\
    \node[2]{B_n}\arrow{sw,l,-}{\langle b\rangle}\arrow{se,l,-}{2^{n-1}}
    \node[2]{P_n}\arrow{sw,r,-}{2^2}\\
    \node{K_n}\arrow{se,r,-}{2^{n-1}}\node[2]{Q_n}\arrow{sw,r,-}{2}\\
    \node[2]{R_n}
  \end{diagram}\]
  \caption{The finite group $G_n$ and its subgroups}\label{fig:Gn}
\end{figure}

\begin{prop}
  These subgroups have the following structure:
  \begin{align*}
    P_n &= (B_{n-1}\times Q_{n-1})\rtimes\langle c,(ac)^4\rangle;\\
    Q_n &= (K_{n-1}\times R_{n-1})\rtimes\langle b,(ac)^4\rangle;\\
    R_n &= (K_{n-1}\times R_{n-1})\rtimes\langle (ac)^4\rangle.
  \end{align*}
\end{prop}
\begin{proof}
  A priori, $P_n$, as a subgroup of $H_n$, maps in $(B_{n-1}\times
  B_{n-1})\rtimes\langle(a,d),(d,a)\rangle$. Restricting to those
  pairs that fix $e_n$ gives the result. Similarly, $Q_n$, as a
  subgroup of $B_n$, maps in $(K_{n-1}\times
  K_{n-1})\rtimes\langle(a,c),(c,a)\rangle$, and $R_n$, as a subgroup
  of $K_n$, maps in $(K_{n-1}\times
  K_{n-1})\rtimes\langle(ac,ca),(ca,ac)\rangle$.
\end{proof}
The group $G_n$ and its subgroups $H_n,B_n,K_n,P_n,R_n,Q_n$ are
arranged in the lattice of Figure~\ref{fig:Gn}, with the quotients or
the indices are represented next to the arrows.

\section{The Group $\tilde G$}\label{sec:Gt}
We describe another fractal group, acting on the same tree $\tree_2$
as $G$. We denote again by $a$ the automorphism permuting the top two
branches, and recursively by $\tilde b$ the automorphism acting as $a$
on the right branch and $\tilde c$ on the left, by $\tilde c$ the
automorphism acting as $1$ on the right branch and $\tilde d$ on the
left, and by $\tilde d$ the automorphism acting as $1$ on the right
branch and $\tilde b$ on the left. In formul\ae,
\begin{alignat*}{2}
  \tilde b(0x\sigma)&=0\overline x\sigma,&\qquad\tilde b(1\sigma)&=1\tilde c(\sigma),\\
  \tilde c(0\sigma)&=0\sigma,&\qquad\tilde c(1\sigma)&=1\tilde d(\sigma),\\
  \tilde d(0\sigma)&=0\sigma,&\qquad\tilde d(1\sigma)&=1\tilde b(\sigma).
\end{alignat*}
Then $\tilde G$ is the group generated by $\{a,\tilde b,\tilde c,\tilde d\}$.
Clearly all these generators are of order $2$, and $\{\tilde b,\tilde
c,\tilde d\}$ is elementary abelian of order $8$.

We write $\tilde H_n=\stab_{\tilde G}(n)$ and $\tilde H=\tilde H_1$.
Explicitly, the map $\phi$ restricts to
\[\phi:\begin{cases}\tilde H\to\tilde H\times\tilde H\\
  \tilde b\mapsto(a,\tilde c),\quad\tilde b^a\mapsto(\tilde c,a)\\
  \tilde c\mapsto(1,\tilde d),\quad\tilde c^a\mapsto(d,1)\\
  \tilde d\mapsto(1,\tilde b),\quad\tilde d^a\mapsto(\tilde b,1).
\end{cases}\]

\begin{prop}
  $\tilde G$ contains elements of finite and infinite order.
\end{prop}
\begin{proof}
  Consider the element $x=a\tilde b\tilde c\tilde d$ of $\tilde G$.
  Then $x^2\in\tilde H$ satisfies $\phi(x^2)=(x,x)$, so
  $x^{2^n}\neq1$ for all $n$; as $\tilde G<\aut(\tree_2)$ has only
  $2$-torsion, it follows that $x$ is of infinite order.
\end{proof}
Note that $x$ acts on $\partial\tree_2$ like an `adding machine'
(see~\cite{bass-:treeaut}). More generally, every spherically
transitive automorphism of $\tree_p$ is conjugated in $\aut(\tree_p)$
to a standard one, called the \emdef{adding machine}, that can be
written $z\mapsto z+1$ after identification of $\partial\tree_p$ with
$\Z_p$.

\begin{prop}
  $\tilde G$ contains $G$ as a subgroup of infinite index.
\end{prop}
\begin{proof}
  The embedding is given by $a\mapsto a$, $b\mapsto\tilde b\tilde d$,
  $c\mapsto\tilde c\tilde b$, $d\mapsto\tilde d\tilde c$. The index is
  infinite because the subgroups $G$ and $\langle a\tilde b\tilde
  c\tilde d\rangle$ do not intersect, one being torsion and the other
  torsion-free.
\end{proof}

Define the elements $u=(a\tilde b)^2$ and $v=(a\tilde d)^2$
in $\tilde G$, and consider its following subgroups:
\[\tilde H=\langle\tilde b,\tilde c,\tilde d\rangle^{\tilde G},\qquad
\tilde B=\langle\tilde b,\tilde d\rangle^{\tilde G},\qquad
\tilde C=\langle\tilde b,v\rangle^{\tilde G},\qquad
\tilde K=\langle u,v\rangle^{\tilde G}.\]

\begin{prop}\label{prop:gtildesub}
  They have the following structure:
  \begin{align*}
    \tilde H &= \langle\tilde b,\tilde c,\tilde d,\tilde b^a,\tilde c^a,\tilde d^a\rangle\text{ is normal of index $2$ in $\tilde G$}.\\
    \tilde B &= \langle\tilde b,\tilde d,\tilde b^a,\tilde d^a,\tilde b^{\tilde ca},\tilde b^{a\tilde ca}\rangle\text{ is normal of index $8$ in $\tilde G$}.\\
    \tilde C &= \langle\tilde b,v,\tilde b^a,\tilde b^{\tilde ca},\tilde b^{a\tilde ca}\rangle\text{ is normal of index $16$ in $\tilde G$}.\\
    \tilde K &= \langle u,v,(a\tilde b\tilde d)^2,u^a,u^{a\tilde c}\rangle\text{ is normal of index $32$ in $\tilde G$}.\\
    \intertext{Furthermore,}
    \phi(\tilde H) &= (\tilde B\times\tilde B)\rtimes\langle\phi(\tilde b),\phi(\tilde b^a)\rangle,\\
    \phi(\tilde B) &= (\tilde C\times\tilde C)\rtimes\langle\phi(\tilde b),\phi(\tilde b^a)\rangle,\\
    \phi(\tilde C) &= (\tilde K\times\tilde K)\rtimes_{\phi\langle([\tilde b,v],[\tilde b^a,v]\rangle}\langle\phi(\tilde b),\phi(\tilde b^a),v\rangle,\\
    \phi(\tilde K) &= (\tilde K\times\tilde K)\rtimes_{\phi\langle[u,v]\rangle^{\tilde G}}\langle u,v\rangle.
  \end{align*}
\end{prop}

\begin{prop}
  $\tilde G$ is spherically transitive, fractal and regular branch
  over its subgroup $\tilde K$.
\end{prop}
\begin{proof}
  $\tilde G$ is fractal by Lemma~\ref{lem:fractal} and the nature of the
  map $\phi$.  As $\tilde K$ is normal, $\phi(\tilde K)$ contains
  $\phi[u,\tilde d]=(1,u)$ and $\phi[u,\tilde c]=(1,v)$, so by
  conjugation it contains $1\times\tilde K$ and $\tilde K\times1$, so
  finally it contains $\tilde K\times\tilde K$.
\end{proof}

\begin{prop}
  $\tilde G$ is just-infinite.
\end{prop}
\begin{proof}
  By direct computation, $[\tilde K:\tilde K']=64$. Apply
  Proposition~\ref{prop:jinf}.
\end{proof}

\begin{prop}
  Define the substitution $\tilde\sigma$ on $\{a,\tilde b,\tilde
  c,\tilde d\}^*$ by
  \[\tilde\sigma:\begin{cases}a\mapsto a\tilde ba,&\tilde b\mapsto\tilde d,\\ \tilde c\mapsto\tilde b,&\tilde d\mapsto\tilde c.\end{cases}\]
  Then $\tilde G$ has a recursive presentation of $L$-type
  \begin{multline}
    \tilde G=\Big\langle a,\tilde b,\tilde c,\tilde d\Big|\,a^2,\tilde
    b^2,\tilde c^2,\tilde d^2,[\tilde b,\tilde c],[\tilde b,\tilde
    d],[\tilde c,\tilde d],\\
    \tilde\sigma^i(a\tilde c)^4,\tilde\sigma^i(a\tilde d)^4,\tilde\sigma^i(a\tilde ca\tilde d)^2,\tilde\sigma^i(a\tilde b)^8,\tilde\sigma^i(a\tilde ba\tilde ba\tilde c)^4,\tilde\sigma^i(a\tilde ba\tilde ba\tilde d)^4,\tilde\sigma^i(a\tilde ba\tilde ba\tilde ca\tilde ba\tilde ba\tilde d)^2\quad(i\ge0)\Big\rangle,
  \end{multline}
  and $\tilde\sigma$ induces an injective expanding endomorphism of
  $\tilde G$ of infinite-index image.
\end{prop}

\begin{proof}
  Consider the groups
  \begin{align*}
    \Gamma &= \langle\alpha,\beta,\gamma,\delta|\,\alpha^2,\beta^2,\gamma^2,\delta^2,[\beta,\gamma],[\beta,\delta],[\gamma,\delta],(\alpha\gamma)^4\rangle,\\
    \Xi &= \langle\beta,\gamma,\delta,\beta^\alpha,\gamma^\alpha,\delta^\alpha\rangle <_2\Gamma.
  \end{align*}
  Then $\tilde G$ is a quotient of $\Gamma$, written $\tilde
  G=\Gamma/\Omega$, via the map $\alpha\mapsto a,\beta\mapsto
  b,\gamma\mapsto c,\delta\mapsto d$, and the map $\phi$ lifts to a
  map $\theta:\Xi\to\Gamma\times\Gamma$. Define
  \[\Omega_n = \{g\in\Gamma|\,\text{$\theta^n$ is applicable and
    }\theta^n(g)=(1,\dots,1)\quad(2^n\text{ copies})\},\]
  where the notation
  implies that $g\in\Xi,\theta(g)\in\Xi\times\Xi,\dots$.  For any word
  $w$ in $\{\alpha,\beta,\gamma,\delta\}$ of length at least $2$
  representing an element of $\Xi$, the corresponding words
  $\theta(w)_{1,2}$ will be strictly shorter; thus every $g\in\Omega$
  eventually gives $1$ through iterated application of $\theta$, and
  thus $\Omega=\cup_{n\ge0}\Omega_n$. We will obtain an explicit
  set of generators for $\Omega_n$: let $\omega_0=(\alpha\gamma)^4$ and
  \[\{\omega_1,\dots,\omega_6\}=\{(\alpha\delta)^4,(\alpha\gamma\alpha\delta)^2,(\alpha\beta)^8,(\alpha\beta\alpha\beta\alpha\gamma)^4,(\alpha\beta\alpha\beta\alpha\delta)^4,(\alpha\beta\alpha\beta\alpha\gamma\alpha\beta\alpha\beta\alpha\delta)^2\}.\]
  Then we claim that for all $n\ge0$
  \[\Omega_n=\langle\tilde\sigma^{j+1}(\omega_0),\tilde\sigma^j(\omega_i)\quad(0\le j\le n-1,1\le i\le6)\rangle^\Gamma.\]
  
  By direct application of the Todd-Coxeter
  algorithm~\cite{gap:manual}, we obtain the presentation
  \[\Xi = \Big\langle\beta,\gamma,\delta,\overline\beta,\overline\gamma,\overline\delta\Big|\,\beta^2,\gamma^2,\delta^2,\overline\beta^2,\overline\gamma^2,\overline\delta^2,[\beta,\gamma],[\beta,\delta],[\gamma,\delta],[\overline\beta,\overline\gamma],[\overline\beta,\overline\delta],[\overline\gamma,\overline\delta],[\gamma,\overline\gamma]\Big\rangle.\]
  
  Computation shows that $\theta(\Xi)$ is of index $8$ in
  $\Gamma\times\Gamma$. From this we obtain, again using Todd-Coxeter,
  the presentation
  \begin{multline*}
    \theta(\Xi)=\Big\langle\beta,\gamma,\delta,\overline\beta,\overline\gamma,\overline\delta\Big|\,
    \beta^2,\gamma^2,\delta^2,\overline\beta^2,\overline\gamma^2,\overline\delta^2,[\beta,\gamma],[\beta,\delta],[\gamma,\delta],[\overline\beta,\overline\gamma],[\overline\beta,\overline\delta],[\overline\gamma,\overline\delta],\\
[\gamma,\overline\gamma],[\gamma,\overline\delta],[\delta,\overline\gamma],[\delta,\overline\delta],(\beta\overline\beta)^4,(\beta\overline\beta x\beta\overline\beta\overline y)^2\,(x,y\in\{\gamma,\delta\})\Big\rangle.
  \end{multline*}
  
  As a consequence, we can write $\ker(\theta)$ as a normal subgroup
  of $\Gamma$ by keeping only those relators of $\theta(\Xi)$ that do
  not appear in $\Xi$ and rewriting them in
  $\{\alpha,\beta,\gamma,\delta\}$, namely
  \[\Omega_1=\ker(\theta)=\langle\omega_1,\dots,\omega_6\rangle^\Gamma.\]

  Then a direct computation shows that
  $\theta\tilde\sigma(\omega_i)=(1,\omega_i)$ for $i=0,\dots,6$. This proves
  that
  \begin{align*}
    \Omega_n&=\{g\in\Xi|\,\theta(g)\in\Omega_{n-1}\times\Omega_{n-1}\}\\
    &=\big(\{\tilde\sigma(g)|\,g\in\Omega_{n-1}\}\cup\Omega_{n-1}\big)^\Gamma\\
    &=\langle\tilde\sigma^{j+1}(\omega_0),\tilde\sigma^j(\omega_i)\quad(0\le j\le n-1,1\le i\le6)\rangle^\Gamma.
  \end{align*}
\end{proof}

\begin{cor}
  All relations of $\tilde G$ have even length. As a consequence, the
  Cayley graph of $\tilde G$ relative to the generating set
  $\{a,\tilde b,\tilde c,\tilde d\}$ is bipartite.
\end{cor}

We believe the relations given in the previous theorem are
independent, and that the method used in~\cite{grigorchuk:bath} can be
used to prove this.

Note that the relations of $G$ can be obtained from those of $\tilde
G$; in the following equalities we indicate by an underscore the
letters affected by a relation in $\tilde G$.
\begin{align*}
  (ad)^4 &= (a\tilde c\tilde d)^4 =_{\tilde G} (a\tilde c\tilde da\tilde d\tilde c)^2 =_{\tilde G} (a\tilde c\tilde da\tilde d(\tilde da)^4\tilde c)(a\tilde c\tilde da\tilde d\tilde c)\\
  &=_{\tilde G} (\underline{a\tilde ca\tilde d}\;\underline{a\tilde da\tilde c})(a\tilde d\tilde ca\tilde c\tilde d) =_{\tilde G} (\tilde da\tilde ca\tilde ca\tilde da)(a\tilde d\tilde ca\tilde c\tilde d) =_{\tilde G} \tilde d(a\tilde c)^4\tilde d =_{\tilde G} \tilde d^2 =_{\tilde G} 1,\\
  \intertext{and}
  (adacac)^4 &= (a\tilde d\underline{\tilde ca\tilde c}\tilde ba\tilde b\tilde c)^4 =_{\tilde G} (\underline{a\tilde da\tilde c}a\tilde ca\tilde ba\tilde b\tilde c)^4 =_{\tilde G} \tilde ca(\tilde d\underline{\tilde ca\tilde ba\tilde ba})(\tilde d\tilde ca\tilde ba\tilde ba)^3a\tilde c\\
  &=_{\tilde G} \tilde ca\big(\tilde d(a\tilde ba\tilde ba\tilde c)^3(\tilde d\tilde ca\tilde ba\tilde ba)^3\big)a\tilde c\\
  &=_{\tilde G} \tilde ca\big(\tilde d(a\tilde ba\tilde ba\tilde c)\underline{(a\tilde ba\tilde ba\tilde c)(a\tilde ba\tilde ba\tilde da\tilde ba\tilde ba)(\tilde c}\tilde da\tilde ba\tilde ba)(\tilde d\tilde ca\tilde ba\tilde ba)\big)a\tilde c\\
  &=_{\tilde G} \tilde ca\big(\tilde d(a\tilde ba\tilde ba\tilde c)\underline{(\tilde da\tilde ba\tilde ba)^2(\tilde d}\tilde ca\tilde ba\tilde ba)\big)a\tilde c =_{\tilde G} \tilde ca\tilde d(a\tilde ba\tilde ba\tilde ca\tilde ba\tilde ba\tilde d)^2\tilde da\tilde c =_{\tilde G} 1.
\end{align*}

\begin{prop}
  The finite quotients $\tilde G_n=\tilde G/\stab_{\tilde G}(n)$ of
  $\tilde G$ have order $2^{13\cdot2^{n-4}+2}$ for $n\ge 4$, and
  $2^{2^n-1}$ for $n\le4$.
\end{prop}
\begin{proof}
  For $n\ge4$, $\phi(\tilde H)$ is a subgroup of index $8$ in $\tilde
  G\times\tilde G$, so $\tilde G_n$ is a subgroup of index $8$ in
  $\tilde G_{n-1}\wr\Z/2$ and $|\tilde G_n|=|\tilde G_{n-1}|^2/4$. For
  $n\le4$ one has $\tilde G_n=\aut(\tree)_n=\Z/2\wr\dots\wr\Z/2$.
\end{proof}

\begin{prop}
  $\tilde K\ge\stab_{\tilde G}(4)$ and $\tilde K'\ge\tilde K'_{(2)}$, so
  $\tilde G$ has the congruence property. Additionally, $\tilde
  K'\ge\stab_{\tilde G}(5)$.
\end{prop}
\begin{proof}
  The first and third assertions can be checked on a computer. For the
  second, $K$ contains $y=[u,d]$ and $z=[u,c]$; these elements satisfy
  $\phi(y)=(1,u)$ and $\phi(z)=(1,v)$. Then $K'$ contains
  $[y,v]=\phi^{-2}(1,1,u,1)$ and $[z,d]=\phi^{-2}(1,1,v,1)$, so it
  contains $\phi^{-2}(1\times1\times K\times1)$ and $K_{(2)}$.
\end{proof}

\begin{cor}
  The closure $\overline{\tilde G}$ of $\tilde G$ in $\aut(\tree)$ is
  isomorphic to the profinite completion $\widehat{\tilde G}$ and is a
  pro-$2$-group. It has Hausdorff dimension $13/16$.
\end{cor}

\subsection{The Growth of $\tilde G$}
By the growth of a group one means the growth, in the sense of
Definition~\ref{defn:growth}, of the group acting on itself.  We
rephrase the definition of growth of a group in a slightly more
general frame:
\begin{defn}
  Let $G$ be a group generated by a finite set $S$, and let
  $\nu:S\to\R_+^*$ be any function. The \emdef{weight} of $g\in G$ is
  \[|g| = \min\{\nu(s_1)+\dots+\nu(s_n)|\,s_1\cdots s_n=g,s_i\in S\}.\]
  The \emdef{growth series} of $G$ with respect to $\nu$ is
  \[F_\nu(\tau) = \sum_{g\in G} e^{\tau|g|}.\]
  This series converges at least in the half-plane
  $\Re(\tau)<-\log(n)/\min_{s\in S}\nu(s)$. Let $\rho(\nu)$, the \emdef{growth
    rate} of $G$ with respect to $\nu$, be the smallest non-positive
  value such that the series converges.
\end{defn}

\begin{prop}
  If $\rho(\nu)<0$, then $G$ has exponential growth, while if
  $\rho(\nu)=0$, then $G$ has intermediate or polynomial growth.
\end{prop}
\begin{proof}
  Let $m$ and $M$ be the minimum and maximum of the weight function
  $\nu$, and set $R=\lim\sqrt[n]{\gamma_G^S(n)}$. By considering
  the series $F_S(\tau)=\sum_{n\ge0}\gamma(n)\tau^n$, whose radius of
  convergence is $1/R$ and comparing it with $F_\nu(\tau)$, we obtain
  \[M\rho(\nu)\le\log(1/R)\le m\rho(\nu),\]
  so $R>1$ is equivalent to $\rho(\nu)<0$.
\end{proof}

The first examples of groups of intermediate growth were constructed
in~\cite{grigorchuk:growth}; the group $G$ is one of them.
\begin{thm}
  $\tilde G$ has intermediate growth.
\end{thm}
\begin{proof}
  First, note that $\tilde G$ cannot have polynomial growth, since it
  contains $G$ whose growth function is greater than $e^{\sqrt
  n}$~\cite{bartholdi:lowerbd}. 

  Take as generators for $\tilde G$ the set $S=\{a,\tilde b,\tilde
  c,\tilde d,\tilde b\tilde c,\tilde b\tilde d,\tilde c\tilde d,\tilde
  b\tilde c\tilde d\}$; let $\theta$ be strictly between the real root
  of the equation $-2+X+X^2+X^3=0$ and $1$, for instance
  $\theta=0.811$ and let $\nu$ be defined by
  \begin{gather*}
    \nu(a)=1,\\
    \nu(\tilde b)=(\theta+\theta^3)/(1-\theta^3)\approx 2.87,\\
    \nu(\tilde c)=(-1+\theta+\theta^2+\theta^3)/(1-\theta^3)\approx 2.14,\\
    \nu(\tilde d)=(\theta^2+\theta^3)/(1-\theta^3)\approx 2.54,\\
    \nu(\tilde b\tilde c)=(-1+\theta+\theta^3)/(1-\theta^3)\approx 0.73,\\
    \nu(\tilde b\tilde d)=(\theta^3)/(1-\theta^3)\approx 1.13,\\
    \nu(\tilde c\tilde d)=(-1+\theta^2+\theta^3)/(1-\theta^3)\approx 0.41,\\
    \nu(\tilde b\tilde c\tilde d)=(1+\theta^3)/(1-\theta^3)\approx 3.28.
  \end{gather*}

  Clearly any element $g\in G$, when expressed as a minimal word in
  $S$, will have the form $[a]x_1ax_2\dots ax_n[a]$, where the first
  and last $a$ are optional and $x_i\in S\setminus\{a\}$. Indeed the
  function $\nu$ satisfies the triangular inequalities
  $\nu(\tilde b)+\nu(\tilde c)<\nu(\tilde b\tilde c)$, etc. Choose
  once and for all a minimal expression for every element of $\tilde G$.
  
  Suppose now for contradiction that $\rho(\nu)<0$.  For some value
  $\eta\in(0,1)$ to be chosen later, partition $\tilde G$ in two sets:
  $A$ containing those elements $g\in G$ whose minimal expression
  $s_1\dots s_n$ contains at least $\eta n$ occurrences of the
  generator $x=\tilde b\tilde c\tilde d$, and $B$ the other elements.
  Define two generating series
  \[F_A(\tau)=\sum_{g\in A}e^{\tau|g|},\qquad F_B(\tau)=\sum_{g\in A}e^{\tau|g|}.\]
  Clearly $F_\nu = F_A+F_B$. We will show that for an appropriate
  value of $\eta$ both $F_A$ and $F_B$ will converge up to some
  $\sigma$ with $\rho(\nu)<\sigma<0$.
  
  We bound $F_A$ by replacing $A$ by a larger set, namely the set of
  all words $s_1\dots s_n$ containing at least $\eta n$ occurrences
  of $x$. Then
  \[F_A(\tau)<\sum_{n\ge0}\binom{n}{\eta n}\left(\sum_{s\in
      S}e^{\tau\nu(s)}\right)^{(1-\eta)n}\left(e^{\tau\nu(x)}\right)^{\eta n}.\]
  By Stirling's formula,
  \[\binom{n}{\eta n}\approx\frac{\sqrt{2\pi\eta(1-\eta)}\sqrt{n}}{\left(\eta^\eta+(1-\eta)^{1-\eta}\right)^n}.\]
  Putting these together, we conclude that $F_A$ converges up to
  any $\sigma>\rho(\nu)$ if
  \[\frac{\left(\sum_{s\in S}e^{\sigma\nu(s)}\right)^{1-\eta}\left(e^{\sigma\nu(x)}\right)^\eta}{\eta^\eta(1-\eta)^{1-\eta}}<1,\]
  and this will hold for $\eta$ large enough, as both the first
  multiplicand and the denominator tend to $1$ as $\eta$ tends to $1$,
  while the second multiplicand tends to $e^{\sigma\nu(x)}<1$.
  
  We then approximate $F_B$ by considering the subset $B'\subset B$ of
  words $s_1\dots s_n$ that either start or end by $a$, but not both;
  and further that contain an even number of $a$s. The series
  $F_{B'}(\tau)$ obtained this way will satisfy $F_B\approx 4F_{B'}$.
  Now $B'$ injects in $G\times G$ through the map $\phi$, written
  $g\mapsto(g_{|0},g_{|1})$. We will compare $|g|$ with
  $|g_{|0}|+|g_{|1}|$. Thanks to the choice of $\nu$, every
  generator $s\neq x$ contributing $\nu(s)$ to $|g|$ will contribute
  at most $\theta\nu(s)$ to $|g_{|0}|+|g_{|1}|$, while every $x$
  contributing $\nu(x)$ to $|g|$ will contribute $\nu(x)$ to
  $|g_{|0}|+|g_{|1}|$. We conclude that for all $g\in B'$ we have
  \[\frac{|g_{|0}|+|g_{|1}|}{|g|}<\frac{\eta\nu(x)+(1-\eta)\min\nu}{\eta\nu(x)+(1-\eta)\theta\min\nu}=:\zeta<1.\]
  This means every element of weight $n$ in $B'$ can be written as a
  pair of elements of $G$ with total weight at most $\zeta n$, or in
  formul\ae
  \[F_{B'}(\tau)\le (F_\nu(\zeta\tau))^2.\]
  The series $F_{B'}$ thus converges up to $\zeta\rho(\nu)>\rho(\nu)$;
  the same holds for $F_B$. Then the series $F_\nu$ converges up to
  $\min(\zeta\rho(\nu),\sigma)>\rho(\nu)$, a contradiction.
\end{proof}

\subsection{The Subgroup $\tilde P$}
Let $e$ be the ray $1^\infty$ and let $\tilde P$ be the corresponding
parabolic subgroup.

\begin{thm}\label{thm:decomPtilde}
  $\tilde P/\tilde P'$ is an infinite elementary $2$-group generated
  by the images of $\tilde b$, $\tilde c=(1,\tilde d)$, $\tilde
  d=(1,\tilde b)$ and of all elements of the form $(1,\dots,1,u^2)$
  or $(1,\dots,1,v)$.  The following decomposition holds:
  \[\tilde P = \bigg(\tilde B\times\Big(\big(\tilde C\times((\tilde K\times\dots)\rtimes\langle u^2,v\rangle)\big)\rtimes\langle\tilde b,u^2,\tilde d,v\rangle\Big)\bigg)\rtimes\langle\tilde b,u^2\rangle.\]
\end{thm}

Define the following subgroups of $\tilde G_n$:
\begin{gather*}
  \tilde B_n = \langle\tilde b,\tilde d\rangle^{\tilde G_n};\qquad
  \tilde C_n = \langle\tilde b,\tilde v\rangle^{\tilde G_n};\qquad
  \tilde K_n = \langle u,v\rangle^{\tilde G_n};\\
  \tilde Q_n = \tilde B_n\cap \tilde P_n;\qquad\tilde R_n = \tilde
  C_n\cap\tilde P_n;\qquad S_n = \tilde K_n\cap \tilde P_n.
\end{gather*}

\begin{figure}
  \[\begin{diagram}
    \node[4]{\tilde G_n}\arrow{s,l,-}{\langle a\rangle}\\
    \node[4]{\tilde H_n}\arrow{sw,l,-}{\langle\tilde c,\tilde c^a\rangle}\arrow{se,l,-}{2^{n-1}}\\
    \node[3]{\tilde B_n}\arrow{sw,l,-}{\langle\tilde d\rangle}\arrow{se,l,-}{2^{n-1}}
    \node[2]{\tilde P_n}\arrow{sw,r,-}{2^2}\\
    \node[2]{\tilde C_n}\arrow{sw,l,-}{\langle\tilde b\rangle}\arrow{se,l,-}{2^{n-1}}
    \node[2]{\tilde Q_n}\arrow{sw,r,-}{2}\\
    \node{K_n}\arrow{se,r,-}{2^{n-1}}\node[2]{R_n}\arrow{sw,r,-}{2}\\
    \node[2]{\tilde S_n}
  \end{diagram}\]
  \caption{The finite group $\tilde G_n$ and its subgroups}\label{fig:Sn}
\end{figure}

\begin{prop}
  These subgroups have the following structure:
  \begin{align*}
    \tilde P_n &= (\tilde B_{n-1}\times\tilde Q_{n-1})\rtimes\langle\tilde b,u^2\rangle;\\
    \tilde Q_n &= (\tilde C_{n-1}\times\tilde R_{n-1})\rtimes\langle\tilde b,u^2\rangle;\\
    \tilde R_n &= (\tilde K_{n-1}\times\tilde S_{n-1})\rtimes_{\langle[b,v]\rangle}\langle b,u^2,v\rangle.
    \tilde S_n &= (\tilde K_{n-1}\times\tilde S_{n-1})\rtimes_{\langle[u^2,v]\rangle}\langle u^2,v\rangle.
  \end{align*}
\end{prop}
\begin{proof}
  The claims match those of Proposition~\ref{prop:gtildesub}, and are
  proved by restricting to elements preserving $e_n$ the `$y$' and
  `$z$' in decompositions of the kind $(x\times y)\rtimes z$.
\end{proof}
The group $\tilde G_n$ and its subgroups $\tilde H_n,\tilde B_n,\tilde
C_n,\tilde K_n,\tilde P_n,\tilde R_n,\tilde Q_n,\tilde S_n$ are
arranged in the lattice of Figure~\ref{fig:Sn}, with the quotients or
the indices are represented next to the arrows.

\section{The Group $\Gamma$}\label{sec:Gamma}
The next three groups we study are subgroups of $\aut(\tree_3)$.
Denote by $a$ the automorphism of $\tree_3$ permuting cyclically the
top three branches. Let $t$ be the automorphisms of $\tree_3$ defined
recursively by
\[t(0x\sigma)=0\overline x\sigma,\qquad t(1x\sigma)=1x\sigma,\qquad
t(2\sigma)=2t(\sigma).\]
Then $\Gamma$ is the subgroup of $\aut(\tree_3)$ generated by
$\{a,t\}$; its growth was studied by Jacek Fabrykowski and Narain
Gupta~\cite{fabrykowski-g:growth2}.

We write $H_n=\stab_\Gamma(n)$ and $H=H_1$.  Explicitly, the map $\phi$
restricts to
\[\phi:\begin{cases}t\to(a,1,t),\quad t^a\to(t,a,1),\quad t^{a^2}\to(1,t,a).\end{cases}\]

Define the elements $x=at$, $y=ta$ of $\Gamma$. Let $K$ be the
subgroup of $\Gamma$ generated by $x$ and $y$, and let $L$ be the
subgroup of $K$ generated by $K'$ and cubes in $K$.

\begin{prop}
  We have the following diagram of normal subgroups:
  
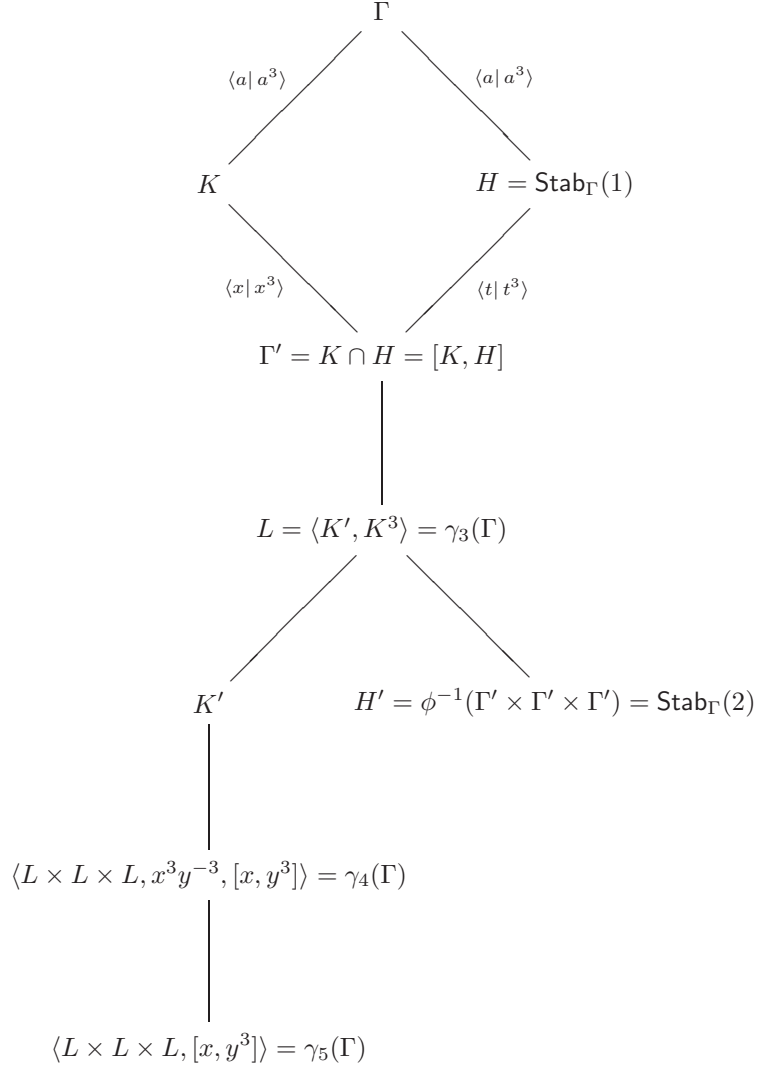
\[\begin{diagram}
    \node[2]{\Gamma}\arrow{sw,l,-}{\langle a|\,a^3\rangle}
    \arrow{se,l,-}{\langle a|\,a^3\rangle}\\
    \node{K}\arrow{se,r,-}{\langle x|\,x^3\rangle}
    \node[2]{H=\stab_\Gamma(1)}\arrow{sw,r,-}{\langle t|\,t^3\rangle}\\
    \smashnode[2]{\Gamma'=K\cap H=[K,H]}\arrow{s,-}\\
    \smashnode[2]{L=\langle K',K^3\rangle=\gamma_3(\Gamma)}\arrow{sw,-}\arrow{se,-}\\
    \node{K'}\arrow{s,-}
    \smashnode[2]{H'=\phi^{-1}(\Gamma'\times\Gamma'\times\Gamma')=\stab_\Gamma(2)}\\
    \smashnode{\langle L\times L\times
    L,x^3y^{-3},[x,y^3]\rangle=\gamma_4(\Gamma)}\arrow{s,-}\\
    \node{\langle L\times L\times
    L,[x,y^3]\rangle=\gamma_5(\Gamma)}
  \end{diagram}\]
  where the quotients are represented next to the arrows; all edges
  represent normal inclusions of index $3$. Furthermore $L=K\cap
  \phi^{-1}(K\times K\times K)$.
\end{prop}
\begin{proof}
  First we prove $K$ is normal in $\Gamma$, of index $3$, by writing
  $y^t=x^{-1}y^{-1}$, $y^{a^{-1}}=y^{-1}x^{-1}$, $y^{t^{-1}}=y^a=x$;
  similar relations hold for conjugates of $x$.  A transversal of $K$
  in $\Gamma$ is $\langle a\rangle$.  All subgroups in the diagram are
  then normal.

  Since $[a,t]=y^{-1}x = t^at^{-1}$, we clearly have $\Gamma'<K\cap H$.
  Now as $\Gamma'\neq K$ and $\Gamma'\neq H$ and $\Gamma'$ has index
  $3^2$, we must have $\Gamma'=K\cap H$. Finally
  $[a,t]=[x,t]^{t^{-1}}$, so $\Gamma'=[K,H]$.

  Next $x^3=[a,t][t,a^{-1}][a^{-1},t^{-1}]$ and similarly for $y$, so
  $K^3<\Gamma'$ and $L<\Gamma'$. Also,
  $\phi[x,y]=(y^{-1},y^{-1},x^{-1})$ and $\phi x^3=(y,x,y)$ both
  belong to $K\times K\times K$, while $[a,t]$ does not; so $L$ is a
  proper subgroup of $\Gamma'$, of index $3$ (since $K/L$ is
  the elementary abelian group $(\Z/3\Z)^2$ on $x$ and $y$).
  
  Consider now $H'$. It is in $\stab_\Gamma(2)$ since $H=\stab_\Gamma(1)$.
  Also, $[t,t^a]=y^3[y^{-1},x]$ and similarly for other conjugates of
  $t$, so $H'<L$, and $\phi[t,t^a]=([a,t],1,1)$, so
  $\phi(H')=\Gamma'\times\Gamma'\times\Gamma'$. Finally $H'$ it is of
  index $3$ in $L$ (since $H/H'=(\Z/3\Z)^3$ on $t,t^a,t^{a^{-1}}$),
  and since $\stab_\Gamma(2)$ is of index $3^4$ in $\Gamma$ (with
  quotient $\Z/3\Z\wr\Z/3\Z$) we have all the claimed equalities.
\end{proof}

\begin{prop}\label{prop:Gammafractal}
  $\Gamma$ is a just-infinite fractal group, is regular branch
  over $\Gamma'$, and has the congruence property.
\end{prop}
\begin{proof}
  $\Gamma$ is fractal by Lemma~\ref{lem:fractal} and the nature of the
  map $\phi$.  By direct computation, $[\Gamma:\Gamma'] =
  [\Gamma':\phi^{-1}(\Gamma'\times\Gamma'\times\Gamma')] =
  [\phi^{-1}(\Gamma'\times\Gamma'\times\Gamma'):\Gamma''] = 3^2$, so
  $\Gamma$ is branched on $\Gamma'$. Then $\Gamma''=\gamma_5(\Gamma)$,
  as is shown in~\cite{bartholdi:lcs}, so $\Gamma''$ has finite index
  and $\Gamma$ is just-infinite by Proposition~\ref{prop:jinf}.

  $\Gamma'\ge\stab_\Gamma(2)$, so $\Gamma$ has the congruence property.
\end{proof}

\begin{prop} We have, with the notation introduced in
  Definition~\ref{defn:subsd},
  \begin{align*}
    \phi(H)&=(\Gamma'\times\Gamma'\times\Gamma')\rtimes_{3-ab}\langle t,t^a,t^{a^2}\rangle,\\
    \phi(\Gamma')&=(\Gamma'\times\Gamma'\times\Gamma')\rtimes_{3-ab}\langle [a,t],[a^2,t]\rangle.
  \end{align*}
\end{prop}

\begin{thm}
  The subgroup $K$ of $\Gamma$ is torsion-free; thus $\Gamma$ is
  virtually torsion-free.
\end{thm}
\begin{proof}
  For $1\neq g\in K$, let $|g|_t$, the $t$-length of $g$, denote the minimal
  number of $t^{\pm1}$'s required to write $g$ as a word over the
  alphabet $\{a^{\pm1},t^{\pm1}\}$. We will show by induction on
  $|g|_t$ that $g$ is of infinite order.

  First, if $|g|_t=1$, i.e.\ $g\in\{x^{\pm1},y^{\pm1}\}$, we
  conclude from $\phi(x^3)=(*,*,x)$ and $\phi(y^3)=(*,*,y)$ that $g$
  is of infinite order.
  
  Suppose now that $|g|_t>1$. If $g\in L$, then
  $\phi(g)=(g_0,g_1,g_2)\in K\times K\times K$ and it suffices to show
  that one of the $g_i$ is of infinite order---this follows by
  induction since $|g_i|_t<|g|_t$ and some $g_i\neq1$. We may thus
  suppose that $g\in K\setminus L$. Up to symmetry, it suffices also
  to consider elements $g$ of the form $\ell x$, $\ell xy$ and $\ell
  xy^{-1}$, for $\ell\in L$. Write
  $\phi(\ell)=(\ell_0,\ell_1,\ell_2)$.

  In the first case, we have $\phi(g^3)=\phi(\ell x)^3=(a\ell_2
  t\ell_1\ell_0,*,*)$. It suffices to show that the first coordinate
  of this expression is non-trivial, as $K$ contains at worst only
  $3$-torsion, being contained in the $3$-Sylow of $\aut(\tree_3)$.
  Now map $\Gamma$ to $\Gamma/\Gamma'$, an elementary abelian group of
  order $9$. One checks that $\ell_0\ell_1\ell_2\equiv1$ in the
  abelian quotient, so the first coordinate maps to $\overline
  a\overline t\not\equiv1$ in $\Gamma/\Gamma'$.
  
  The second case is handled in the same way. Finally, if $g=\ell
  xy^{-1}$, then $\phi(g^3)\in L\times L\times L$, so $\phi^2(g^3)\in
  K\times\dots\times K$ ($9$ factors); each factor has strictly
  smaller $t$-length than $g$, and as before the projection in one of
  the coordinates onto the abelian quotient gives some $x\not\equiv1$.
\end{proof}

\begin{prop}
  The finite quotients $\Gamma_n=\Gamma/H_n$ of $\Gamma$ have order
  $3^{3^{n-1}+1}$ for $n\ge2$, and $3$ for $n=1$.
\end{prop}
\begin{proof}
  Follows immediately from $[\Gamma:\Gamma']=3^2$ and
  $[\Gamma':\phi^{-1}(\Gamma'\times\Gamma'\times\Gamma')]=3^2$.
\end{proof}

\begin{cor}
  The closure $\overline\Gamma$ of $\Gamma$ in $\aut(\tree)$ is
  isomorphic to the profinite completion $\widehat\Gamma$ and is a
  pro-$3$-group. It has Hausdorff dimension $1/3$.
\end{cor}

\subsection{The Subgroup $P$}
Let $e$ be the infinite sequence $2^\infty$, and let $P$ be the 
corresponding parabolic subgroup.

\begin{thm}\label{thm:decomPG}
  $P/P'$ is an infinite elementary $3$-group generated by $t$,
  $t^a$ and all elements of the form $(1,\dots,1,[a,t])$. The
  following decomposition holds:
  \[P = \bigg(\Gamma'\times\Gamma'\times \Big(\big(\Gamma'\times\Gamma'\times ((\Gamma'\times\Gamma'\times\dots)\rtimes_{3-ab}\langle[a,t]\rangle)\big)\rtimes_{3-ab}\langle[a,t]\rangle\Big)\bigg)\rtimes_{3-ab}\langle t,t^a\rangle,\]
  where each factor (of nesting $n$) in the decomposition acts on the
  subtree just below some $e_n$ but not containing $e_{n+1}$.
\end{thm}

Define the following subgroups of $\Gamma_n$:
\[\Gamma'_n = \langle[a,t]\rangle^{\Gamma_n};\qquad Q_n = \Gamma'_n\cap P_n.\]

\begin{prop}
  These subgroups have the following structure:
  \begin{align*}
    P_n &= (\Gamma'_{n-1}\times\Gamma'_{n-1}\times Q_{n-1})\rtimes_{3-ab}\langle t,t^a\rangle;\\
    Q_n &= (\Gamma'_{n-1}\times\Gamma'_{n-1}\times Q_{n-1})\rtimes_{3-ab}\langle[a,t]\rangle.
  \end{align*}
\end{prop}

\section{The Group $\overline\Gamma$}\label{sec:Gammab}
Recall $a$ denotes the automorphism of $\tree_3$ permuting
cyclically the top three branches. Let now $t$ be the
automorphism of $\tree_3$ defined recursively by
\[t(0x\sigma)=0\overline x\sigma,\qquad t(1x\sigma)=1\overline
x\sigma,\qquad t(2\sigma)=2t(\sigma).\]
Then $\overline\Gamma$ is the subgroup of $\aut(\tree_3)$ generated by
$\{a,t\}$.

We write $H_n=\stab_{\overline\Gamma}(n)$ and $H=H_1$. Explicitly, the
map $\phi$ restricts to
\[\phi:\begin{cases}t\to(a,a,t),\quad t^a\to(t,a,a),\quad t^{a^2}\to(a,t,a).\end{cases}\]

Define the elements $x=ta^{-1}$, $y=a^{-1}t$ of $\overline\Gamma$, and let $K$
be the subgroup of $\overline\Gamma$ generated by $x$ and $y$. Then $K$ is
normal in $\overline\Gamma$, because $x^t=y^{-1}x^{-1}$, $x^a=x^{-1}y^{-1}$,
$x^{t^{-1}}=x^{a^{-1}}=y$, and similar relations hold for conjugates of $y$.
Moreover $K$ is of index $3$ in $\overline\Gamma$, with transversal $\langle a\rangle$.

\begin{lem}\label{lem:GammaHK}
  $H$ and $K$ are normal subgroups of index $3$ in $\overline\Gamma$, and
  $\overline\Gamma'=\stab_K(1)=H\cap K$ is of index $9$; furthermore $\phi(H\cap
  K)\triangleleft K\times K\times K$. For any element
  $g=(u,v,w)\in\phi(H\cap K)$ one has $wvu\in H\cap K$.
\end{lem}  
\begin{proof}
  First note that $\stab_K(1)=\langle x^3,y^3,xy^{-1},y^{-1}x\rangle$,
  for every word in $x$ and $y$ whose number of $a$'s is divisible by
  $3$ can be written in these generators. Then compute
  \begin{alignat*}{2}
    \phi(x^3)&=(y,x^{-1}y^{-1},x),&\qquad\phi(y^3)&=(x^{-1}y^{-1},x,y),\\
    \phi(xy^{-1})&=(1,x^{-1},x),&\qquad\phi(y^{-1}x)&=(y,1,y^{-1}).
  \end{alignat*}
  The last assertion is also checked on this computation.
\end{proof}

\begin{prop}
  Writing $c=[a,t]=x^{-1}y^{-1}x^{-1}$ and $d=[x,y]$, we have the
  following diagram of normal subgroups:
  
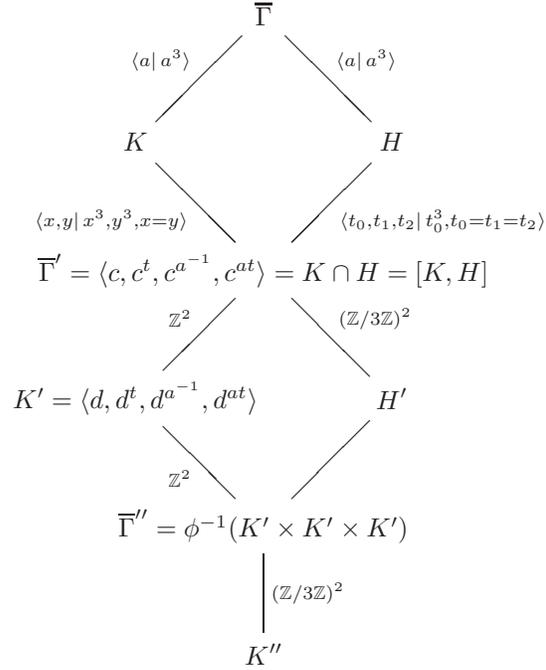
\[\begin{diagram}
    \node[2]{\overline\Gamma}\arrow{sw,l,-}{\langle a|\,a^3\rangle}
    \arrow{se,l,-}{\langle a|\,a^3\rangle}\\
    \node{K}\arrow{se,r,-}{\langle x,y|\,x^3,y^3,x=y\rangle}
    \node[2]{H}\arrow{sw,r,-}{\langle t_0,t_1,t_2|\,t_0^3,t_0=t_1=t_2\rangle}\\
    \smashnode[2]{\overline\Gamma'=\langle c, c^t, c^{a^{-1}}, c^{at} \rangle=K\cap H=[K,H]}\arrow{sw,t,-}{\Z^2}\arrow{se,t,-}{(\Z/3\Z)^2}\\
    \smashnode{K'=\langle d, d^t, d^{a^{-1}}, d^{at} \rangle}\arrow{se,b,-}{\Z^2}
    \node[2]{H'}\arrow{sw,-}\\
    \smashnode[2]{\overline\Gamma''=\phi^{-1}(K'\times K'\times
      K')}\arrow{s,r,-}{(\Z/3\Z)^2}\\
    \node[2]{K''}
  \end{diagram}\]
  where the quotients are represented next to the arrows; additionally,
  \begin{gather*}
    K/K' = \langle x,y|\,[x,y]\rangle\cong\Z^2,\\
    \overline\Gamma'/\overline\Gamma'' = \langle
    c,c^t,c^{a^{-1}},c^{at}|\,[c,c^t],\dots\rangle\cong \Z^4,\\
    K'/K'' = \langle d,d^t,d^{a^{-1}},d^{at}|\,[d,d^t],\dots,(d/d^{at})^3,(d^{a^{-1}}/d^t)^3\rangle\cong\Z^2\times(\Z/3\Z)^2.
  \end{gather*}

  Writing each subgroup in the generators of the groups above it, we
  have
  \begin{gather*}
    K = \langle x=at^{-1},y=a^{-1}t\rangle,\\
    H = \langle t,t_1=t^a,t_2=t^{a^{-1}}\rangle,\\
    \overline\Gamma' = \langle b_1=xy^{-1},b_2=y^{-1}x,b_3=x^3,b_4=y^3\rangle
    = \langle c_1=tt_1^{-1},c_2=tt_1t,c_3=tt_2^{-1},c_4=tt_2t\rangle.
  \end{gather*}
\end{prop}

\begin{prop}\label{prop:GammaBfractal}
  $\overline\Gamma$ is a fractal group, is weakly branch, and
  just-nonsolvable; however it is not branch.
\end{prop}
\begin{proof}
  $\overline\Gamma$ is fractal by Lemma~\ref{lem:fractal} and the
  nature of the map $\phi$. The subgroup $K$ described above has an
  infinite-index derived subgroup $K'$ (with infinite cyclic
  quotient), from which we conclude that $\overline\Gamma$ is not
  just-infinite; indeed $K'$ is normal
  in $\overline\Gamma$ and
  $\overline\Gamma/K'\cong\Z^2\rtimes\left(\begin{smallmatrix}-1&1\\
      -1&0\end{smallmatrix}\right)$ is infinite.
\end{proof}

\begin{prop}
  The subgroup $K$ of $\overline\Gamma$ is torsion-free; thus
  $\overline\Gamma$ is virtually torsion-free.
\end{prop}
\begin{proof}
  For $1\neq g\in K$, let $|g|_t$, the $t$-length of $g$, denote the minimal
  number of $t^{\pm1}$'s required to write $g$ as a word over the
  alphabet $\{a^{\pm1},t^{\pm1}\}$. We will show by induction on
  $|g|_t$ that $g$ is of infinite order.

  First, if $|g|_t=1$, i.e.\ $g\in\{x^{\pm1},y^{\pm1}\}$, we
  conclude from $\phi(x^3)=(*,*,x)$ and $\phi(y^3)=(*,*,y)$ that $g$
  is of infinite order.
  
  Suppose now that $|g|_t>1$, and $g\in H_n\setminus H_{n+1}$. Then
  there is some sequence $\sigma$ of length $n$ that is fixed by $g$
  and such that $g_{|\sigma}\not\in H$. By Lemma~\ref{lem:GammaHK},
  $g_{|\sigma}\in K$, so it suffices to show that all $g\in K\setminus
  H$ are of infinite order.
  
  Such a $g$ can be written as $\phi^{-1}(u,v,w)z$ for some
  $(u,v,w)\in \phi(K\cap H)$ and $z\in\{x^{\pm1},y^{\pm1}\}$; by
  symmetry let us suppose $z=x$. Then
  $g^3=\phi^{-1}(uavawt,vawtua,wtuava)=\phi^{-1}(g_0,g_1,g_2)$, say.
  For any $i$, we have $|g_i|_t\le|g|_t$, because all the components
  of $\phi(x)$ and $\phi(y)$ have $t$-length $\le1$. We distinguish
  three cases:
  \begin{enumerate}
  \item $g_i=1$ for some $i$. Then consider the image $\overline{g_i}$
    of $g_i$ in $\overline\Gamma/\overline\Gamma'$. By
    Lemma~\ref{lem:GammaHK}, $wvu\in G'$, so
    $\overline{g_i}=1=\overline{a^2t}$. But this is a contradiction,
    because $\overline\Gamma/\overline\Gamma'$ is elementary abelian
    of order $9$, generated by the independent images $\overline a$
    and $\overline t$.
  \item $0<|g_i|_t<|g|_t$ for some $i$. Then by induction
    $g_i$ is of infinite order, so $g^3$ too, and $g$ too.
  \item $|g_i|_t=|g|_t$ for all $i$. We repeat the
    argument with $g_i$ substituted for $g$. As there are finitely many
    elements $h$ with $|h|_t=|g|_t$, we will eventually
    reach either an element of shorter length or an element
    already considered. In the latter case we obtain a relation
    of the form $\phi^n(g^{3^n})=(\dots,g,\dots)$ from which $g$ is
    seen to be of infinite order.
  \end{enumerate}
\end{proof}

\begin{prop}
  The finite quotients $\overline\Gamma_n=\overline\Gamma/H_n$ of
  $\overline\Gamma$ have order $3^{\frac14(3^n+2n+3)}$ for $n\ge2$,
  and $3^{\frac12(3^n-1)}$ for $n\le2$.
\end{prop}
\begin{proof}
  Define the following family of two-generated finite abelian groups:
  \[A_n = \begin{cases}\langle x,y|\,x^{3^{n/2}},y^{3^{n/2}},[x,y]\rangle&\text{ if }n\equiv0[2],\\
    \langle x,y|\,x^{3^{(n+1)/2}},y^{3^{(n+1)/2}},(xy^{-1})^{3^{(n-1)/2}},[x,y]\rangle&\text{ if }n\equiv1[2].
  \end{cases}\]
  First suppose $n\ge 2$; Consider the diagram of groups described
  above, and quotient all the groups by $H_n$. Then the quotient
  $K/K'$ is isomorphic to $A_n$, generated by $x$ and $y$, and the
  quotient $K'/\overline\Gamma''$ is isomorphic to $A_{n-1}$,
  generated by $[x,y]$ and $[x,y]^t$. As $|A_n|=3^n$, the index of
  $K_n'$ in $\overline\Gamma_n$ is $3^{n+1}$ and the index of
  $\overline\Gamma_n''$ is $3^{2n}$.  Then as
  $\overline\Gamma''_n\cong K_{n-1}^3$ and $|\overline\Gamma_2''|=1$
  we deduce by induction that
  $|\overline\Gamma''_n|=3^{\frac14(3^n-6n+3)}$ and
  $|K'_n|=3^{\frac14(3^n-2n-1)}$, from which
  $|\overline\Gamma_n|=3^{2n}+|\overline\Gamma''_n|=3^{\frac14(3^n+2n+3)}$
  follows.
  
  For $n\le2$ we have $\overline\Gamma_n=\aut(\tree)_n=\Z/3\wr\dots\wr\Z/3$.
\end{proof}

\begin{cor}
  The closure $\overline{\,\overline\Gamma\,}$ of $\overline\Gamma$ in
  $\aut(\tree)$ has Hausdorff dimension $1/2$.
\end{cor}

\begin{prop}
  We have
  \begin{align*}
    \phi(H) &= (K'\times K'\times K')\rtimes_A\langle
    t_0,t_1,t_2\rangle,\\
    \phi(K') &= (K'\times K'\times
    K')\rtimes_B\langle
    d,d^t\rangle,
  \end{align*}
  where $A$ is such that $\langle
  t_0,t_1,t_2\rangle/A\cong\Z^4\rtimes\Z/3\Z$ and $B$ is such that
  $\langle d,d^t\rangle/B\cong\Z^2$.
\end{prop}

\subsection{The Subgroup $P$}
Let $e$ be the infinite sequence $2^\infty$, and let $P$ be the 
corresponding parabolic subgroup.

\begin{thm}\label{thm:decomPGB}
  $P/P'$ is the direct product of $(\Z/3\Z)^2$ (generated by $t$ and
  $atat^{-1}a$) and an infinitely-generated free abelian group,
  generated by $[tt_1t,tt_2t]$. The following decomposition holds:
  \[P = \bigg(K'\times K'\times \Big(\big(K'\times K'\times ((K'\times
  K'\times\dots)\rtimes\langle[tt_1t,tt_2t]\rangle)\big)\rtimes\langle[tt_1t,tt_2t]\rangle\Big)\bigg)\rtimes\langle t,t_1t_2^{-1}\rangle,\]
  where each factor (of nesting $n$) in the decomposition acts on the
  subtree just below some $e_n$ but not containing $e_{n+1}$.
\end{thm}

Define the following subgroups of $\Gamma_n$:
\[K'_n = {\langle x,y\rangle'}^{\Gamma_n};\qquad Q_n = K'_n\cap P_n.\]

\begin{prop}
  These subgroups have the following structure:
  \begin{align*}
    P_n &= (K'_{n-1}\times K'_{n-1}\times Q_{n-1})\rtimes_{3-ab}\langle (\Z^4\ltimes\Z/3\Z)\rangle;\\
    Q_n &= (K'_{n-1}\times K'_{n-1}\times Q_{n-1})\rtimes\Z^2\rangle.
  \end{align*}
\end{prop}

\section{The Group $\doverline\Gamma$}\label{sec:Gammabb}
Recall $a$ denotes the automorphism of $\tree_3$ permuting
cyclically the top three branches. Let now $t$ be the
automorphism of $\tree_3$ defined recursively by
\[t(0x\sigma)=0\overline x\sigma,\qquad
t(1x\sigma)=1\doverline x\sigma,\qquad
t(2\sigma)=2t(\sigma).\]

Then $\doverline\Gamma$ is the subgroup of $\aut(\tree_3)$
generated by $\{a,t\}$; it was studied by Narain Gupta and Said
Sidki~\cite{gupta-s:burnside,gupta-s:infinitep,sidki:subgroups,sidki:pres}.

We will use the following known facts:
\begin{thm}
  $\doverline\Gamma$ is a torsion $3$-group.
\end{thm}

\begin{prop}
  We have the following diagram of normal subgroups:
  \[\begin{diagram}
    \node{\doverline\Gamma}\arrow{s,r,-}{\langle a|\,a^3\rangle}\\
    \node{H=\stab_{\doverline\Gamma}(1)}\arrow{s,r,-}{\langle t|\,t^3\rangle}\\
    \node{\doverline\Gamma'=[\doverline\Gamma,H]}\arrow{s,r,-}{[a,t]}\\
    \node{\gamma_3(\doverline\Gamma)=\doverline\Gamma^3=\stab_{\doverline\Gamma}(2)}\arrow{s,r,-}{(at)^3}\\
    \node{H'=\phi^{-1}(\doverline\Gamma'\times\doverline\Gamma'\times\doverline\Gamma')}
  \end{diagram}\]
  where the quotients are represented next to the arrows; all edges
  represent normal inclusions of index $3$.
\end{prop}
\begin{proof}
  Clearly $H$ is normal of index $3$, being the kernel of the
  epimorphism $a\to a,t\to1$. Then $\doverline\Gamma'\neq H$ (as can
  be checked in the finite quotient $\doverline\Gamma_2$) but is of
  index at most $3^2$, so has precisely that index. Moreover,
  $\doverline\Gamma'$ is generated by the $[a^{\pm1},t^{\pm1}]$: one
  has $[a,t]^a=[a^{-1},t][a,t]^{-1}$, $[a,t]^t=[a,t]^{-1}[a,t^{-1}]$,
  etc.
  
  $\gamma_3(\doverline\Gamma)<\doverline\Gamma^3$ holds in all
  $3$-groups, and $\doverline\Gamma^3$ has index $3^3$ because it is
  $2$-generated $2$-step nilpotent.
  
  Now consider $H'$. It is in $\stab_{\doverline\Gamma}(2)$ since
  $H=\stab_{\doverline\Gamma}(1)$.  Also,
  $[t,t^a]=(ta)^3(a^{-1}ta^{-1})^3$ and similarly for other
  conjugates, so $H'<\doverline\Gamma^3$, and
  $\phi[t^{-a^2}t^{-a},t^{-a}t^{-1}]=([a,t],1,1)$, so
  $\phi(H')=\doverline\Gamma'\times\doverline\Gamma'\times\doverline\Gamma'$.
  Finally $H'$ it is of index $3$ in $\doverline\Gamma^3$ (since
  $H/H'=(\Z/3\Z)^3$ on $t,t^a,t^{a^{-1}}$), and since
  $\stab_{\doverline\Gamma}(2)$ is of index $3^4$ in $\Gamma$ (with
  quotient $\Z/3\Z\wr\Z/3\Z$) we have all the claimed equalities.
\end{proof}

\begin{prop}\label{prop:GammaBBfractal}
  $\doverline\Gamma$ is a just-infinite fractal group, and is a
  regular branch group over $\doverline\Gamma'$.
\end{prop}
\begin{proof}
  $\doverline\Gamma$ is fractal by Lemma~\ref{lem:fractal} and the
  nature of the map $\phi$.  By direct computation,
  $[\doverline\Gamma:\doverline\Gamma'] =
  [\doverline\Gamma':\phi^{-1}(\doverline\Gamma'\times\doverline\Gamma'\times\doverline\Gamma')]=[\phi^{-1}(\doverline\Gamma'\times\doverline\Gamma'\times\doverline\Gamma'):\doverline\Gamma'']=3^2$,
  so $\doverline\Gamma$ is branched on $\doverline\Gamma'$ and is
  just-infinite by Proposition~\ref{prop:jinf}.
\end{proof}

\begin{prop}
  $\doverline\Gamma'\ge\stab_{\doverline\Gamma}(2)$, so
  $\doverline\Gamma$ has the congruence property.
\end{prop}

\begin{prop} We have
  \begin{align*}
    \phi(H)&=(\doverline\Gamma'\times\doverline\Gamma'\times\doverline\Gamma')\rtimes_{3-ab}\langle t,t^a,t^{a^2}\rangle,\\
    \phi(\doverline\Gamma')&=(\doverline\Gamma'\times\doverline\Gamma'\times\doverline\Gamma')\rtimes_{3-ab}\langle [a,t],[a^2,t]\rangle.
  \end{align*}
\end{prop}

\subsection{The Subgroup $P$}
Let $e$ be the infinite sequence $2^\infty$, and let $P$ be the
corresponding parabolic subgroup.

\begin{thm}\label{thm:decomPGBB}
  $P/P'$ is an infinite elementary $3$-group generated by $t$,
  $t^at^{a^2}$ and all elements of the form $(1,\dots,1,tt^at^{a^2})$. The
  following decomposition holds:
  \[P = \bigg(\doverline\Gamma'\times\doverline\Gamma'\times
  \Big(\big(\doverline\Gamma'\times\doverline\Gamma'\times
  ((\doverline\Gamma'\times\doverline\Gamma'\times\dots)\rtimes_{3-ab}\langle
  tt^at^{a^2}\rangle)\big)\rtimes_{3-ab}\langle tt^at^{a^2}\rangle\Big)\bigg)\rtimes_{3-ab}\langle t,t^at^{a^2}\rangle,\]
  where each factor (of nesting $n$) in the decomposition acts on the
  subtree just below some $e_n$ but not containing $e_{n+1}$.
\end{thm}

Define the following subgroups of $\doverline\Gamma_n$:
\[\doverline\Gamma'_n = \langle[a,t]\rangle^{\doverline\Gamma_n};\qquad Q_n = \doverline\Gamma'_n\cap P_n.\]

\begin{prop}
  These subgroups have the following structure:
  \begin{align*}
    P_n &= (\doverline\Gamma'_{n-1}\times\doverline\Gamma'_{n-1}\times Q_{n-1})\rtimes_{3-ab}\langle t,t^at^{a^2}\rangle;\\
    Q_n &= (\doverline\Gamma'_{n-1}\times\doverline\Gamma'_{n-1}\times
    Q_{n-1})\rtimes_{3-ab}\langle tt^at^{a^2}\rangle.
  \end{align*}
\end{prop}

\section{Quasi-Regular Representations}\label{sec:qr}
In this section we show how the information we gathered on the groups
and their subgroups yields results on their representations.  For $G$
a group acting on a tree and $P$ its parabolic subgroup, we let
$\rho_{G/P}$ denote the quasi-regular representation of $G$ on the
space $\ell^2(G/P)$.

First of all consider the infinite-dimensional representations
$\rho_{G/P}$. The criterion of irreducibility for quasi-regular
representations was discovered by George Mackey and is as follows (the
definition of commensurator is given after the theorem's statement):
\begin{thm}[Mackey~\cite{mackey:representations,burger-h:irreducible}]\label{thm:mackey}
  Let $G$ be an infinite group and let $P$ be any subgroup of $G$.
  Then the quasi-regular representation $\rho_{G/P}$ is irreducible if
  and only if $\comm_G(P)=P$.
\end{thm}

\begin{defn}
  The \emdef{commensurator} (also called \emph{quasi-normalizer}) of a
  subgroup $H$ of $G$ is
  \[\comm_G(H) = \{g\in G|\,H\cap H^g\text{ is of finite index in
    }H\text{ and }H^g\}.\]
  Equivalently, letting $H$ act on the left on the right cosets $\{gH\}$,
  \[\comm_G(H) = \{g\in G|\,H\cdot(gH)\text{ and }H\cdot(g^{-1}H)\text{ are finite orbits}\}.\]
\end{defn}
The equivalence follows, for $T$ a finite transversal, from
\[H=\bigsqcup_{t\in T\subset H}t\cdot(H\cap H^g)\Longleftrightarrow
HgH=\bigsqcup_{t\in T\subset H}t\cdot gH.\]

\begin{prop}\label{prop:commP}
  If $G$ is weakly branch, then $\comm_G(P)=P$.
\end{prop}
\begin{proof}
  Take $g\in G\setminus P$, with $P=\stab_G(e)$ for some ray $e$; we
  will show that $P\cap P^g$ is of infinite index in $P^g$. Let $n$ be
  such that $\sigma:=e_1\dots e_n\neq g(e_1\dots e_n)$. Then
  $\rist_{P^g}(\sigma)=\rist_G(\sigma)$, while by
  Lemma~\ref{lem:inforbits} the index of $\rist_{P\cap
    P^g}(\sigma)=\rist_P(\sigma)$ in $\rist_G(\sigma)$ is infinite.
\end{proof}

\begin{cor}\label{cor:irreducible}
  If $G$ is weakly branch, then $\rho_{G/P}$ is irreducible.
\end{cor}

The quasi-regular representations we consider are good approximants of
the regular representation in the following sense:
\begin{thm}
  $\rho_G$ is a subrepresentation of $\bigotimes_{P\text{ parabolic
  }}\rho_{G/P}$.
\end{thm}
\begin{proof}
  Since $\bigcap_{g\in G}P^g=1$, it follows that the $G$-space $G$ is
  a subspace of $\prod_{g\in G}G/P_g$. The representation on a product
  of spaces is the tensor product of the representation on the spaces.
\end{proof}

We have a continuum of parabolic subgroups $P_e=\stab_G(e)$, where $e$
runs through the boundary of a tree, so formally we also have a
continuum of quasi-regular representations.
If $G$ is countable, there are uncountably many non-equivalent
representations, because among the uncountably many $P_e$ only
countably many are conjugate.
We therefore have the
\begin{thm}
  There are uncountably many non-equivalent representations of the
  form $\rho_{G/P}$, where $P$ is a parabolic subgroup.
\end{thm}

We now consider the finite-dimensional representations $\rho_{G/P_n}$,
where $P_n$ is the stabilizer of the vertex at level $n$ in the ray
defining $P$. These are permutational representations on the sets
$G/P_n$. The $\rho_{G/P_n}$ are factors of the representation
$\rho_{G/P}$.  Noting that $P=\bigcap_{n\ge0}P_n$, it follows that
\[\rho_{G/P_n}\Rightarrow\rho_{G/P},\]
in the sense that for any non-trivial $g\in G$ there is an $n\in\N$
with $\rho_{G/P_n}(g)\neq1$.

\subsection{Hecke Algebras}
Corollary~\ref{cor:irreducible} showed that the quasi-regular
representation $\rho_{G/P}$ is irreducible for all of our examples. We
now describe the decomposition of the finite quasi-regular
representations $\rho_{G/P_n}$. It turns out that it is closely
related to the orbit structure of $P_n$ on $G/P_n$, through the
\emph{Hecke algebra}. The result we shall prove is:
\begin{thm}\label{thm:finitedec}
  $\rho_{G/P_n}$ and $\rho_{\tilde G/\tilde P_n}$ decompose as a direct
  sum of $n+1$ irreducible components, one of degree $2^i$ for each
  $i\in\{1,\dots,n-1\}$ and two of degree $1$.
  
  $\rho_{\Gamma/P}$, $\rho_{\overline\Gamma/P}$ and
  $\rho_{\doverline\Gamma/P}$ decompose as a direct sum of $2n+1$
  irreducible components, two of degree $2^i$ for each
  $i\in\{1,\dots,n-1\}$ and three of degree $1$.
\end{thm}
The proof of this theorem will appear after the following definitions
and lemmata.

\begin{defn}
  Let $G$ be a group and $P$ a subgroup. Set $Q=\comm_G(P)$, and
  define
  \[\C[G,P]=\Big\{f:Q\to\C\Big|\,f(pqp')=f(q)\forall p,p'\in
  P\text{ and }\operatorname{supp}(f)\subset\bigcup_{\text{finite}}PqP\Big\},\]
  i.e. those $(P,P)$-invariant functions on $Q$ whose support is
  a finite union of $(P,P)$-double cosets. $\C[G,P]$ is an algebra for
  the convolution product
  \[(f\cdot g)(x) = \sum_{y\in G/P} f(xy) g(y^{-1}).\]
  
  The \emdef{Hecke algebra} (also called the \emdef{intersection
    algebra}) $\Hecke(G,P)$ is the weak closure of $\C[G,P]$ in
  $\mathcal L(\ell^2(G/P))$.
\end{defn}
A few remarks are in order. First, the convolution product is well
defined on $\C[G,P]$, since every double coset $PqP$ is a finite union
of left (or right) cosets. Second, $\Hecke(G,P)$ coincides with
the commutant $\rho_{G/P}(G)'$ of the right-regular representation of
$G$ in $\mathcal L(\ell^2(G/P))$. That $\Hecke(G,P)$ commutes with
$\rho'$ is obvious, since these two operators derive from left- and
right-actions on $G$. That $\Hecke(G,P)$ is the full commutant
requires an argument, based on approximation of functions in
$\Hecke(G,P)$ by finite-support functions.

Third, the whole theory of Hecke algebra can be extended to locally
compact $G$ and compact-open $P$ --- see for
instance~\cite{tzanev:phd}. One then defines $\C[G,P]$ as those
bi-$P$-invariant continuous maps $G\to\C$ whose support is contained
in a finite union of $PJP$, where the $J$ are compact-open subgroups
of $G$. This algebra is represented in $\mathcal L(L^2(G/P),\mu)$,
where $\mu$ is the projection of the Haar measure to $G/P$ (which,
beware, need not be $G$-invariant!). We shall not make use of this
theory.

A variant of this notion, which we will use, is obtained by taking $G$
profinite and $P$ closed. Then $\C[G,P]$ consists of those
bi-$P$-invariant continuous maps $G\to\C$ whose support is contained
in a finite union of $PJP$, where the $J$ are neighbourhoods of the
identity in $G$.

$\Hecke(G,P)$ is topologically spanned by compactly supported
$(P-P)$-biinvariant functions on $G$. The following result stresses
the importance of the Hecke algebra in the study of representation
decomposition: \def\0{\cite[Section~11D]{curtis-r:methods}}
\begin{thm}[\0]\label{thm:hecke}
  Suppose $[G:P]$ is finite. Then $\Hecke(G,P)$ is a semi-simple
  algebra.  There is a canonical bijection between isotypical
  components of $\rho_{G/P}$ and simple factors of $\Hecke(G,P)$,
  which maps $\chi^n$ (for $\chi$ simple) to $M_n(\C)$.
\end{thm}
Then, if $\Hecke(G,P)$ is abelian, its decomposition in simple modules
is has as many components as there are double cosets $PgP$ in $G$.

In our examples, the spaces have the following order of magnitude: the
core of $P_n$ is the normal subgroup $H_n=\bigcap_{g\in G} P_n^g$, of
index $\sim e^{e^n}$. The subgroup $P_n$ is of index $\sim e^n$. The
number of double cosets is $\sim n$. We give the precise results for
our five examples.

\subsection{Orbits In $G/P_n$}\label{subs:computedcosets}
As the double cosets $P_ngP_n$ are in one-to-one correspondence with
the orbits of $P_n$ on $G/P_n$ we shall now describe the orbits for
this action.

\begin{lem}\label{lem:Korbit}
  There are two $K_n$-orbits on $\Sigma^n$: those sequences starting
  with $0$ and those starting with $1$.

  $P_n$ has $n+1$ orbits in $\Sigma^n$; they are $1^n$ and the
  $1^i0\Sigma^{n-1-i}$ for $0\le i<n$. The orbits of $P$ in
  $\tree_\Sigma$ are the $1^i0\Sigma^*$ for all $i\in\N$.
\end{lem}
\begin{proof}
  As $K_n$ contains $K_{n-1}\times K_{n-1}$, it follows by induction
  that $K_n$ acts transitively on the sets $00\Sigma^{n-2}$ and
  $01\Sigma^{n-2}$. As $K_n$ contains $(ab)^2=(ca,ac)$, it also permutes
  $00\Sigma^{n-2}$ and $01\Sigma^{n-2}$, so it acts transitively on
  $0\Sigma^{n-1}$. The same holds for $1\Sigma^{n-1}$.

  The last assertion follows from Theorem~\ref{thm:decomP}.
\end{proof}

\begin{lem}
  There are two $\tilde K_n$-orbits on $\Sigma^n$: those sequences
  starting with $0$ and those starting with $1$.

  $\tilde P_n$ has $n+1$ orbits in $\Sigma^n$; they are $1^n$ and the
  $1^i0\Sigma^{n-1-i}$ for $0\le i<n$. The orbits of $\tilde P$ in
  $\tree_\Sigma$ are the $1^i0\Sigma^*$ for all $i\in\N$.
\end{lem}
\begin{proof}
  Completely similar to~\ref{lem:Korbit}.
\end{proof}

\begin{lem}\label{lem:Gporbit}
  There are three $\Gamma'_n$-orbits on $\Sigma^n$: those sequences starting
  with $0$, those starting with $1$ and those starting with $2$.

  $P_n$ has $2n+1$ orbits in $\Sigma^n$; they are $2^n$ and the
  $2^i0\Sigma^{n-1-i}$ and $2^i1\Sigma^{n-1-i}$ for $0\le i<n$. The
  orbits of $P$ in $\tree_\Sigma$ are the $2^i0\Sigma^*$ and
  $2^i1\Sigma^*$ for all $i\in\N$.
\end{lem}
\begin{proof}
  As $\Gamma'_n$ contains
  $\Gamma'_{n-1}\times\Gamma'_{n-1}\times\Gamma'_{n-1}$, it follows by
  induction that $\Gamma'_n$ acts transitively on the sets
  $00\Sigma^{n-2}$, $01\Sigma^{n-2}$ and $02\Sigma^{n-2}$. As
  $\Gamma'_n$ contains $[a,t]=(ta^{-1},a,t^{-1})$, it also permutes
  $00\Sigma^{n-2}$, $01\Sigma^{n-2}$ and $02\Sigma^{n-2}$, so it acts
  transitively on $0\Sigma^{n-1}$. The same holds for $1\Sigma^{n-1}$
  and $2\Sigma^{n-1}$.

  The last assertion follows from Theorem~\ref{thm:decomPG}
\end{proof}

\begin{lem}\label{lem:GBporbit}
  For the group $\overline\Gamma$, there are three $K'_n$-orbits on
  $\Sigma^n$: those sequences starting with $0$, those starting with
  $1$ and those starting with $2$.

  $P_n$ has $2n+1$ orbits in $\Sigma^n$; they are $2^n$ and the
  $2^i0\Sigma^{n-1-i}$ and $2^i1\Sigma^{n-1-i}$ for $0\le i<n$. The
  orbits of $P$ in $\tree_\Sigma$ are the $2^i0\Sigma^*$ and
  $2^i1\Sigma^*$ for all $i\in\N$.
\end{lem}
\begin{proof}
  As $K'_n$ contains $K'_{n-1}\times K'_{n-1}\times K'_{n-1}$, it
  follows by induction that $K'_n$ acts transitively on the sets
  $00\Sigma^{n-2}$, $01\Sigma^{n-2}$ and $02\Sigma^{n-2}$. As $K'_n$
  contains $[x,y]=(at,at,ta)$, it also permutes $00\Sigma^{n-2}$,
  $01\Sigma^{n-2}$ and $02\Sigma^{n-2}$, so it acts transitively on
  $0\Sigma^{n-1}$. The same holds for $1\Sigma^{n-1}$ and
  $2\Sigma^{n-1}$.

  The last assertion follows from Theorem~\ref{thm:decomPGB}
\end{proof}

\begin{lem}\label{lem:GBBporbit}
  There are three $\doverline\Gamma'_n$-orbits on $\Sigma^n$: those
  sequences starting with $0$, those starting with $1$ and those
  starting with $2$.

  $P_n$ has $2n+1$ orbits in $\Sigma^n$; they are $2^n$ and the
  $2^i0\Sigma^{n-1-i}$ and $2^i1\Sigma^{n-1-i}$ for $0\le i<n$. The
  orbits of $P$ in $\tree_\Sigma$ are the $2^i0\Sigma^*$ and
  $2^i1\Sigma^*$ for all $i\in\N$.
\end{lem}
\begin{proof}
  Completely similar to~\ref{lem:Gporbit}.
\end{proof}

\subsection{Gelfand Pairs}
We have seen the Hecke algebra $\Hecke(G,P_n)$ is roughly of
dimension $n$. Its structure is further simplified by the following
consideration:
\begin{defn}[\cite{diaconis:representations}]
  Let $G$ be a group and $P$ any subgroup. The pair $(G,P)$ is a
  \emdef{Gelfand pair} if all irreducible subrepresentations of
  $\rho_{G/P}$ have multiplicity $1$.
\end{defn}

\def\0{\cite[Exercise~18, page~306]{curtis-r:methods}}
\def\1{\cite[Theorem~1.20]{mackey:representations}}
\begin{lem}[\0,\1]
  $(G,P)$ is a Gelfand pair if and only if $\Hecke(G,P)$ is abelian.
\end{lem}

\begin{prop}\label{prop:gelfand}
  In our five examples the pairs $(G,P_n)$ form a Gelfand pair for all
  $n\in\N$.
\end{prop}
\begin{proof}
  Clearly $P_0=G$ so $\Hecke(G,P_0)=\C$ is abelian. Furthermore,
  $P_{n+1}$ is a subgroup of $P_n$, and the natural map
  $G/H_{n+1}\twoheadrightarrow G/H_n$ induces a map
  $P_{n+1}/H_{n+1}\twoheadrightarrow P_n/H_n$, so
  $\Hecke(G,P_n)\cong\Hecke(G/H_n,P_n/H_n)$ is a direct summand of
  $\Hecke(G,P_{n+1})$, and their dimensions differ by $d-1$, which is
  $1$ or $2$ (recall $d$ is the degree of the regular tree on which
  $G$ acts). Now writing
  \[\Hecke(G,P_{n+1}) = \Hecke(G,P_n) \oplus A,\]
  we see that $A$ is semi-simple and of dimension $d-1<4$. All such
  semisimple algebras are abelian, $A\cong\C^{d-1}$, so
  $\Hecke(G,P_{n+1})$ is abelian too.
\end{proof}

\begin{proof}[Proof of Theorem~\ref{thm:finitedec}]
  By Proposition~\ref{prop:gelfand}, the Hecke algebra $\Hecke(G,P_n)$
  is abelian, so it is isomorphic to $\C^{N_n}$, where $N_n$ is its
  dimension. This $N_n$ in turn is equal to the number of double
  cosets $P_ngP_n$.  These numbers $N_n$ are computed in the
  corollaries in Subsection~\ref{subs:computedcosets}. By
  Theorem~\ref{thm:hecke}, the number of irreducible
  subrepresentations of $\rho_{G/P_n}$ is $N_n$. Finally,
  $\rho_{G/P_n}=\rho_{G/P_{n-1}}\oplus A_{n,1}\oplus\dots\oplus
  A_{n,d-1}$, where the $A_{n,i}$ are irreducible representations.
  Since $\dim\rho_{G/P_n}=d^n$ and $\dim A_i$ is a power of $d$, the
  only possibility is that $\dim A_{n,i}=d^{n-1}$ for all
  $i\in\{1,\dots,d-1\}$, and
  \[\rho_{G/P_n} = \rho_{G/P_0}\oplus A_{1,1}\oplus\dots\oplus
  A_{1,d-1}\oplus\dots\oplus A_{n,1}\oplus\dots\oplus
  A_{n,d-1}.\]
\end{proof}

It may well be that for all $GGS$ groups the Hecke algebra associated
to a parabolic subgroup is commutative.


\section{Acknowledgments}
The authors are immensely grateful to Professor Pierre de la Harpe
who, by inviting the second author for a trimester in Geneva,
facilitated the work in which the results presented here were
obtained, and to Professor Marc Burger, who invited the second author
for a stay in the ETH in Z\"urich during which this paper was
completed.

\bibliography{mrabbrev,people,math,grigorchuk,bartholdi}
\end{document}